\documentclass[11pt]{article}

\begin{document}

\newtheorem{condition}{Condition}
\newcommand{\bt}{\begin{theorem}}
\newcommand{\et}{\end{theorem}}
\newcommand{\beaa}{\begin{eqnarray*}}
\newcommand{\eeaa}{\end{eqnarray*}}
\newcommand{\eps}{\epsilon}
\newcommand{\al}{\alpha}
\newcommand{\bet}{\beta}
\newcommand{\kap}{\kappa}
\newcommand{\Del}{\Delta}

\newcommand{\calA}{{\cal A}}
\newcommand{\calB}{{\cal B}}
\newcommand{\calC}{{\cal C}}
\newcommand{\calD}{{\cal D}}
\newcommand{\calF}{{\cal F}}
\newcommand{\calG}{{\cal G}}
\newcommand{\calH}{{\cal H}}
\newcommand{\calL}{{\cal L}}
\newcommand{\calM}{{\cal M}}
\newcommand{\calP}{{\cal P}}
\newcommand{\calR}{{\cal R}}
\newcommand{\calS}{{\cal S}}
\newcommand{\calT}{{\cal T}}
\newcommand{\calU}{{\cal U}}
\newcommand{\calV}{{\cal V}}
\newcommand{\calX}{{\cal X}}
\newcommand{\calY}{{\cal Y}}
\newcommand{\clb}{{\cal B}}
\newcommand{\clp}{{\cal P}}
\newcommand{\cle}{{\cal E}}
\newcommand{\clc}{{\cal C}}
\newcommand{\cll}{{\cal L}}

\newcommand{\til}{\tilde}
\newcommand{\essinf}{{\rm ess\, inf}}
\newcommand{\w}{\wedge}
\newcommand{\lt}{\left}
\newcommand{\rt}{\right}
\newcommand{\pl}{\partial}
\newcommand{\id}{{\rm id}}

\def\theequation{\arabic{section}.\arabic{equation}}
\def\thetheorem{\arabic{section}.\arabic{theorem}}
\def\theass{\arabic{section}.\arabic{ass}}
\def\thecond{\arabic{section}.\arabic{cond}}

\newcounter{bean}
\newcommand{\benuma}{\setlength{\labelwidth}{.25in}
\begin{list}%
{(\alph{bean})}{\usecounter{bean}}}
\newcommand{\eenuma}{\end{list}}

\newcommand{\beginsec}{\setcounter{equation}{0}}

% square for the end of proof
\newcommand{\ink}{\rule{.5\baselineskip}{.55\baselineskip}}

\newtheorem{theorem}{Theorem}[section]
\newtheorem{remark}[theorem]{Remark}
\newtheorem{lemma}[theorem]{Lemma}
\newtheorem{cor}[theorem]{Corollary}
\newtheorem{defn}[theorem]{Definition}
\newtheorem{cond}[theorem]{Condition}
\newtheorem{assu}[theorem]{Condition}
\newtheorem{prop}[theorem]{Proposition}
\newtheorem{pr}[theorem]{Property}
\newtheorem{definition}[theorem]{Definition}

\newcommand{\noi}{\noindent }
\newcommand{\und}{\underline }
\newcommand{\rimply}{\Rightarrow }
\newcommand{\la}{\lambda  }
\newcommand{\conv}{{\rm conv}}
\newcommand{\cone}{{\rm cone}}
\newcommand{\ti}{\tilde }
\newcommand{\cc}{\cdot }
\newcommand{\non}{\nonumber }
\newcommand{\dist}{{\rm dist}}
\newcommand{\clf}{{\cal F}}
\newcommand{\ii}{\iota}
\newcommand{\In}{{\rm In}}
\newcommand{\pa}{\partial}
\newcommand{\ba}{\begin{array}}
\newcommand{\ea}{\end{array}}
\newcommand{\bea}{\begin{eqnarray}}
\newcommand{\eea}{\end{eqnarray}}
\newcommand{\beas}{\begin{eqnarray*}}
\newcommand{\eeas}{\end{eqnarray*}}
\newcommand{\be}{\begin{equation}}
\newcommand{\ee}{\end{equation}}
\newcommand{\bc}{\begin{center}}
\newcommand{\ec}{\end{center}}
\newcommand{\ben}{\begin{enumerate}}
\newcommand{\een}{\end{enumerate}}
\newcommand{\lan}{\langle}
\newcommand{\ran}{\rangle}
\newcommand{\ei}{\end{itemize}}

\newcommand{\prf}{{\bf Proof.\ }}
\newcommand{\ds}{\displaystyle}
\newcommand{\overbar}{\bar}
\newcommand{\ve}{\varepsilon}
\newcommand{\N}{I\!\!N}
\newcommand{\Z}{Z\!\!Z}
\newcommand{\goto}{\rightarrow}

\newcommand{\brmk}{\begin{remark}\begin{em}}
\newcommand{\ermk}{\end{em}\end{remark}}
\newcommand{\bexa}{\begin{example}\per\begin{em}}
\newcommand{\eexa}{\end{em}\end{example}}

\newcommand{\skp}{\vspace{\baselineskip}}
\newcommand{\Pf}{{\bf Proof.\ }}
\newcommand{\E}{I\!\!E}
\newcommand{\IP}{I\!\!P}
\newcommand{\scal}{{\cal {S}}}
\newcommand{\acal}{{\cal {A}}}
\newcommand{\NN}{{\cal {N}}}
\newcommand{\I}{{\cal {I}}}
\newcommand{\PP}{{\cal {P}}}
\newcommand{\F}{{\cal {F}}}
\newcommand{\C}{{\cal {C}}}
\newcommand{\G}{{\cal {G}}}
\newcommand{\D}{\Delta}
\newcommand{\B}{{\cal {B}}}
\newcommand{\A}{{\cal {A}}}
\newcommand{\Om}{\Omega}
\newcommand{\om}{\omega}
\newcommand{\sig}{\sigma}
\newcommand{\del}{\delta}
\newcommand{\Lsc}{{\cal L}}
\newcommand{\R}{I\! \! R}
\newcommand{\INT}{I\! \! N_0}
\newcommand{\ov}{\overline}
\newcommand{\Df}{\doteq}
\newcommand{\rarr}{\rightarrow}
\newcommand{\larr}{\leftarrow}
\newcommand{\add}{\addtocounter{J}{1}}
\newcommand{\half}{\frac{1}{2}}
\newcommand{\beq}{\begin{eqnarray*}}
\newcommand{\eeq}{\end{eqnarray*}}
\newcommand{\beqn}{\begin{eqnarray}}
\newcommand{\eeqn}{\end{eqnarray}}
\newcommand{\inn}[2]{{#1}\cdot {#2}}
\newcommand{\INN}[2]{\left \langle {#1}, {#2}\right \rangle}
\newcommand{\M}{{\cal M}}
\newcommand{\Su}{{\cal S}}
\newcommand{\ep}{\epsilon}

%%%%%%%%%%%%%%%%%%%%%%%%%%%%%%%%%%% start of logo
{
\centerline{\hbox{
\vrule height -0.4 pt depth 0.8 pt width 26.5 em
\kern - 26.5 em
\raise  0.03ex  \hbox{\bf E} 
\raise  0.06ex \hbox{l} 
\raise .13ex \hbox{e}
\raise .24ex \hbox{c} 
\raise .45ex \hbox{t} 
\raise .78ex \hbox{r} 
\raise 1.31ex \hbox{o} 
\raise 2.08ex \hbox{n} 
\raise 3.14ex \hbox{i} 
\raise 4.53ex \hbox{c} 
\kern   1em
\raise 8.15ex \hbox{\bf J} 
\raise 10.15ex \hbox{o} 
\raise 12.04ex \hbox{u} 
\raise 13.60ex \hbox{r} 
\raise 14.64ex \hbox{n} 
\kern .3 em
\vrule
\kern -.3em
\raise 15ex \hbox{a} 
\raise 14.64ex \hbox{l} 
\kern   1em
\raise 12.04ex \hbox{o} 
\raise 10.15ex \hbox{f} 
\kern   1em
\raise 6.23ex \hbox{\bf P} 
\raise 4.53ex \hbox{r} 
\raise 3.14ex \hbox{o} 
\raise 2.08ex \hbox{b} 
\raise 1.31ex \hbox{a} 
\raise .78ex \hbox{b} 
\raise .45ex \hbox{i} 
\raise .24ex \hbox{l} 
\raise .13ex \hbox{i} 
\raise .06ex \hbox{t} 
\raise .03ex \hbox{y}
 }}}
%%%%%%%%%%%%%%%%%%%%%%%%%%%%%%%%%%%%%%% end of logo

\bigskip
\centerline{Vol. 7 (2002) Paper no. 12, pages 1--36.}
\bigskip
\centerline{Journal URL}
\centerline{http://www.math.washington.edu/\~{}ejpecp/}
\centerline{Paper URL}
\centerline 
{http://www.math.washington.edu/\~{}ejpecp/EjpVol7/paper12.abs.html}

\vskip0.5truein

\centerline{\bf STABILITY PROPERTIES OF CONSTRAINED JUMP-DIFFUSION PROCESSES}

\bigskip

\begin{center}

{\bf Rami Atar}\\
Department of Electrical Engineering, Technion, Haifa 32000, Israel\\
%{\tt atar@ee.technion.ac.il}\\
\bigskip
{\bf Amarjit Budhiraja}\\
Department of Statistics, University of North Carolina,
Chapel Hill, NC 27599-3260, USA\\
%{\tt budhiraj@email.unc.edu}

\end{center}

\bigskip

\noindent{\bf Abstract}:
We consider a class of jump-diffusion processes,
constrained to a polyhedral cone $G\subset\R^n$, where
the constraint vector field is constant on each face of the
boundary.
The constraining mechanism corrects for ``attempts'' of
the process to jump outside the domain.
Under Lipschitz continuity of the Skorohod map $\Gamma$,
it is known that there is a cone $\calC$ such that
the image $\Gamma\phi$ of a deterministic linear
trajectory $\phi$ remains bounded if and only if $\dot\phi\in\calC$.
Denoting the generator of a corresponding unconstrained jump-diffusion
by $\cll$, we show that a key condition for the process
to admit an invariant probability measure is that for $x\in G$,
$\cll\,\id(x)$ belongs to a compact subset of $\calC^o$.

Stability cone. Harris recurrence.

\bigskip

\bigskip

\noindent
Submitted to EJP on September 24, 2001. Final version accepted on 
March 20, 2002.

\noindent
This research was supported in part by the US-Israel
Binational Science Foundation and the fund for the
promotion of research at the Technion (RA),
and the IBM junior faculty development award, University of North
Carolina at Chapel Hill (AB).

\thispagestyle{empty}

\newpage

\section{Introduction}\label{sec:intro}

In this work we consider stability properties of a class of
jump-diffusion processes that are constrained to lie in a convex
closed polyhedral cone.
Let $G$ be a cone in $\R^n$, given as the intersection $\cap_i G_i$
of the half spaces
$$ G_i = \{ x \in \R^n :\inn{x}{n_i} \ge 0 \}, \quad i=1,\ldots, N,$$
where $n_i$, $i=1,\ldots, N$ are given unit vectors.
It is assumed that the origin is a proper vertex of $G$, in the
sense that there exists a closed half space $G_0$ with
$G\cap G_0=\{0\}$. Equivalently, there exists a unit vector
$a_0$ such that
\be\label{asn:1}
\{x\in G:x\cdot a_0\le1\}
\ee
is compact. Note that, in particular, $N\ge n$.
Let $F_i=\pl G\cap\pl G_i$. With each
face $F_i$ we associate a unit vector $d_i$ (such that
$\inn{d_i}{n_i} > 0$). This vector defines the {\it direction of
constraint} associated with the face $F_i$.
The {\em constraint vector field} $d(x)$ is defined for $x\in\pl G$ as
the set of all unit vectors in the cone generated by
$\{d_i, i\in\In(x)\}$,
where
$$\In(x) \Df\{i \in \{1,\ldots, N\}: \inn{x}{n_i} = 0 \}.
$$
Under further assumptions on $(n_i)$ and $(d_i)$, one can define
a {\em Skorohod map} $\Gamma$ in the space of right continuous
paths with left limits, in a way which is consistent with the
constraint vector field $d$. Namely, $\Gamma$ maps a path $\psi$
to a path $\phi=\psi+\eta$ taking values in $G$, so that $\eta$
is of bounded variation, and, denoting the total variation of $\eta$ on
$[0,s]$ by $|\eta|(s)$, $d\eta(\cdot)/d|\eta|(\cdot)\in
d(\phi(\cdot))$. The precise definition of $\Gamma$ and the
conditions assumed are given in Section \ref{sec:results}. The
constrained jump-diffusion studied in this paper is the second
component $Z$ of the pair $(X,Z)$ of processes satisfying \bea
\label{eq:X}
  X_t=z_0+\int_0^t\beta(Z_s)ds+\int_0^t a(Z_s)dW_s
  +\int_{[0,t]\times E}h(\del(Z_{s-},z))[N(ds,dz)-q(ds,dz)]
  \nonumber \\
  +\int_{[0,t] \times E}h'(\del(Z_{s-},z))N(ds,dz),
\eea \be \label{eq:Z}
  Z=\Gamma(X).
\ee
Here, $W$ and $N$ are the driving $m$-dimensional Brownian
motion and Poisson random measure on $\R_+\times E$; $\beta$, $a$
and $\del$ are (state-dependent)
coefficients and $h$ is a truncation function (see
Section \ref{sec:results} for definitions and assumptions). For
illustration, consider as a special case of (\ref{eq:X}),
(\ref{eq:Z}), the case where $X$ is a L\'evy process with
piecewise constant paths and finitely many jumps over finite time
intervals. Then $X_t=x+\sum_{s\le t}\Del X_s$, where $\Del
X_s=X_s-X_{s-}$. In this case, $Z$ is given as $Z_t=x+\sum_{s\le
t}\Del Z_s$, where $\Del Z_s$ can be defined recursively in a
straightforward way. Namely, if $Z_{s-}+\Del X_s\in G$, then
$\Del Z_s=\Del X_s$. Otherwise, $Z_s=Z_{s-}+\Del X_s+\al d$,
where $\al\in(0,\infty)$, $Z_s\in\pl G$, and $d\in d(Z_s)$. In
general, this set of conditions may not have a solution $(\al,d)$,
or may
have multiple solutions. However, the assumptions we put on the
map $\Gamma$ will ensure that this
recursion is uniquely solvable, and as a result, that
the process $Z$ is well defined.

A related model for which recurrence and transience properties
have been studied extensively is that of a semimartingale
reflecting Brownian motion (SRBM) in polyhedral cones \cite{BD,
dupwil, harwil, harwil2, HoRo}. Roughly speaking, a SRBM is a
constrained version, using a ``constraining mechanism'' as
described above, of a Brownian motion with a drift. In a recent
work \cite{abd}, sufficient conditions for positive recurrence of
a constrained diffusion process with a state dependent drift and
(uniformly nondegenerate) diffusion coefficients were obtained.
Under the assumption of regularity of the map $\Gamma$ (as in
Condition \ref{regular} below), it was shown that if the drift
vector field takes values in the cone $\clc$ generated by the
vectors $-d_i$, $i=1,\ldots, N$, and stays away, uniformly, from
the boundary of the cone, then the corresponding constrained
diffusion process is positive recurrent and admits a unique
invariant measure.  The technique used there critically relies on
certain estimates on the exponential moments of the constrained
process.  The current work aims at showing that $\calC$ plays the
role of a stability cone in a much more general setting of
constrained jump-diffusions for which only the first moment is
assumed to be finite. The natural definition of the drift vector
field in the case of a jump-diffusion is $\ti \beta \Df
\cll\,\id$, where $\cll$ denotes the generator of a related
``unconstrained'' jump-diffusion (see (\ref{neti})), and $\id$
denotes the identity mapping on $\R^n$. In the case of a L\'evy
process with finite mean, the drift is simply
$\til\beta(x)=E_xX_1-x$ (which is independent of $x$). Our basic
stability assumption is that the range of $\til\beta$ is contained
in $\cup_{k\in\N}\,k\calC_1$, where $\calC_1$ is a compact subset
of the interior of $\calC$. Under this assumption, our main
stability result states (Theorem \ref{th:rec}): There exists a
compact set $A$ such that for any compact $C\subset G$, \be
\label{intro1} \sup_{x \in C} E_x\tau_A < \infty , \ee where
$\tau_A$ is the first time $Z$ hits $A$, and $E_x$ denotes the
expectation under which $Z$ starts from $x$.  The proof of this
result is based on the construction of a Lyapunov function, and on
a careful separate analysis of small and large jumps of the Markov
process.  As another consequence of the existence of a Lyapunov
function we show that $Z$ is bounded in probability. From the
Feller property of the process it then follows that it admits at
least one invariant measure.  Finally, under further suitable
communicability conditions (see Conditions \ref{commune} and
\ref{resolve}) it follows that the Markov process is positive
Harris recurrent and admits a unique invariant measure.

The study of these processes is motivated by problems in
stochastic network theory (see \cite{kushheavy} for a review). The
assumptions we make on the Skorohod map are known to be satisfied
by a large class of applications, including single class open
queueing networks (see \cite{dupish1}, \cite{harrei}).

For a sampling of stability results on constrained processes with
jumps we list \cite{cheman, dai, KaKe, Ke, mal, MeDo}.
We take an approach similar to that of
\cite{dupwil}, where the stability
properties of SRBM in an orthant are proved by means of
constructing a Lyapunov function.
At the cost of putting conditions that guarantee strong
existence and uniqueness of solutions to the SDE, we are able
to treat diffusions with jumps and state-dependent coefficients.
One of the key properties of the Lyapunov function $f$ constructed
in \cite{dupwil}, is that $Df(x)\cdot b \le -c<0$ for $x \in
G\setminus\{0\}$, where $b$ denotes the constant drift vector of
the unconstrained driving Brownian motion.  In a state-dependent
setting, an analogous condition must hold simultaneously for all
$b$ in the range of $\til\beta$. The construction of the Lyapunov
function is therefore much more involved. The basic stability
assumption referred to above plays a key role in this
construction.

The paper is organized as follows.  In Section 2 we present basic
definitions, assumptions, statements of main results and their
corollaries.  Section 3 is devoted to the proof of (\ref{intro1}),
under the assumption that a suitable Lyapunov function exists.  We
also show in this section that the Markov process is bounded in
probability.  In Section 4 we present the construction of the Lyapunov
function.  Since many arguments are similar to those in \cite{dupwil},
we have tried to avoid repetition wherever possible.  Finally, we have
included certain standard arguments in the appendix for the sake of
completeness.

The following notation is used in this paper.
The boundary relative to $\R^n$
of a set $A \subset \R^n$ is denoted by $\partial
A$. The convex hull of $A$ is denoted by $\conv(A)$.
The cone $\{\sum_{i\in I}\al_iv_i:\al_i\ge0, i\in I\}$ generated
by $(v_i,i\in I)$, is denoted by $\cone\{v_i, i\in I\}$.
The open ball of radius $r$ about $x$
is denoted by $B(x,r)$, and the unit sphere in $\R^n$ by $S^{n-1}$.
$D([0,\infty):\R^n)$ denotes the space of functions mapping
$[0,\infty)$ to $\R^n$ that are right continuous and have limits
from the left.  We endow $D([0,\infty):\R^n)$ with the usual
Skorohod topology. We define
$D_A([0,\infty):\R^n) \doteq \{ \psi \in
D([0,\infty):\R^n) : \psi(0) \in A\}. $ For $\eta \in
D([0,\infty):\R^n)$, $| \eta | (T)$ denotes the total
variation of $\eta$ on $\left[ 0,T \right]$ with respect to the
Euclidean norm on $\R^n$. The Borel $\sig$-field on
$\R^n$ is denoted by $\clb(\R^n)$ and the space of
probability measures on $(\R^n, \clb(\R^n))$ by
$\clp(\R^n)$. Finally, $\al$ denotes a positive constant, whose value
is unimportant and may change from line to line.

\section{Setting and results}\label{sec:results}
\beginsec
Recall from Section \ref{sec:intro} the assumptions on the set $G$
and the definition of the vector field $d$,
$$ d(x) \Df \cone\{d_i, i\in\In(x)\}\cap S^{n-1}.
$$
For $x \in \pa G$, define the set
$n(x)$ of inward normals to $G$ at $x$ by $$ n(x) \Df \{\nu: |\nu| =
1, \;\;\; \nu\cdot(x-y) \le 0, \;\; \forall \; y \in G\}.$$
Let $\Lambda$
be the collection of all the subsets of $\{1, 2, \ldots, N\}$. We
will make the following basic assumption regarding the vectors
$(d_i, n_i)$.
\begin{assu}
\label{compS} For each $\lambda \in \Lambda$, $\lambda \neq
\emptyset$, there exists a vector $d^{\lambda}\in\cone\{d_i,i\in\lambda\}$
with
\be \inn{d^{\lambda}}{n_i} > 0 \;\;\; \mbox{for
all}\;\; i \in \lambda. \ee
\end{assu}
\begin{remark}\rm
\label{dual} An important consequence (cf. \cite{dupwil}) of the
above assumption is that for each $\lambda \in \Lambda$,$\lambda
\neq \emptyset$ there exists a vector $n^{\lambda}$ such that
$n^{\lambda} \in n(x)$ for all $x \in G$ satisfying $\In(x) =
\lambda$ and \be \inn{n^{\lambda}}{d_i} > 0 \;\;\; \mbox{for
all}\;\; i \in \lambda. \ee
\end{remark}

\begin{defn}
\label{def-sp} Let $\psi \in D_G([0,\infty):\R^n)$ be given. Then
$(\phi, \eta) \in D([0,\infty):\R^n)\times D([0,\infty):\R^n)$
 solves the Skorohod problem (SP) for $\psi$
with respect to $G$ and $d$  if and only if $\phi(0) = \psi (0)$,
 and for all
$t \in [0, \infty)$ (1) $\phi(t) = \psi (t) + \eta (t)$; (2) $\phi
(t) \in G$; (3) $| \eta | (t) < \infty$; (4) $\ds | \eta | (t) =
\int_{[0,t]}
 I_{\left\{ \phi (s) \in \partial G \right\} }
d | \eta | (s)$; (5) There exists (Borel) measurable $\gamma : [
0, \infty) \rightarrow \R^{k}$ such that $\gamma (t) \in d( \phi
(t))$ ($d|\eta |$-almost everywhere) and
$ \eta (t) = \int_{[0,t]}
\gamma (s) d | \eta | (s). $
\end{defn}
On the domain $D \subset  D_G([0,\infty):\R^n)$ on which there is
a unique solution to the SP we define the Skorohod
map (SM) $\Gamma$ as $\Gamma(\psi) \Df \phi$, if $(\phi, \psi -
\phi)$
 is the unique
solution of the SP posed by $\psi$. We will make
the following assumption on the regularity of the SM
defined by the data $\{(d_i, n_i); i = 1, 2, \ldots, N\}$.
\begin{assu}
\label{regular} The SM is well defined on all of
$D_G([0,\infty):\R^n)$, i.e., $D = D_G([0,\infty):\R^n)$ and the
SM is Lipschitz continuous in the following sense. There exists a
constant $\ell < \infty$ such that for all $\phi_1, \phi_2 \in
D_G([0,\infty):\R^n)$: \be \sup_{0 \le t <
\infty}|\Gamma(\phi_1)(t) - \Gamma(\phi_2)(t)| \le \ell \sup_{0 \le t <
\infty}|\phi_1(t) - \phi_2(t)|. \label{lip} \ee
\end{assu}
We will assume without loss of generality that $\ell \ge 1$. We
refer the reader to \cite{dupish1,dupram1,harrei}
for sufficient conditions for this regularity property to hold.

We now introduce the constrained processes that will be studied in
this paper.
\begin{defn}
\label{reflevy} Let $(X_t)$ be a L\'{e}vy process starting from
zero (i.e. $X_0 = 0$) , with the L\'{e}vy measure $K$ on $(\R^n,
\clb(\R^n))$. Define a ``constrained L\'{e}vy process'', starting
from $z_0 \in G$, by the relation
$$
Z\Df \Gamma(z_0 + X).
$$
\end{defn}
Recall that a L\'{e}vy measure $K$ is a measure that
satisfies the condition
$\int_{\R^n} |y|^2 \wedge 1 K(dy) < \infty$ (see \cite{Ber},
Chapter 1). We will make one
additional assumption on $K$, as
follows.
\begin{assu}
\label{finmean} The L\'{e}vy measure $K$ satisfies
$$
\int_{\R^n} |y| 1_{|y| \ge 1} K(dy) < \infty .$$
\end{assu}
The above assumption holds if and only if
the L\'{e}vy process $X_t$
 has finite
mean.

We now define the reflected jump-diffusions considered in this
work. On a complete filtered probability space
$(\Om,\calF,(\calF_t),P)$, let an $m$-dimensional standard
Brownian motion $W$ and a Poisson random measure $N$ on
$\R_+\times E$, with intensity measure $q(dt,dz)=dt\otimes F(dz)$
be given.  Here, $(E, \cle)$ is a Blackwell space
and $F$ is a positive $\sig$-finite measure on
$(E,\cle)$. For all practical purposes, $(E, \cle)$ can be taken to be
$(\R^n, \clb(\R^n))$ (see \cite{jacshi}).  Let a truncation function
$h:\R^n\to\R^n$ be a continuous bounded function
satisfying $h(x)=x$ is a neighborhood of the origin and with
compact support. We fix such a function throughout, and
denote also $h'(x) \Df x-h(x)$.
The reflected
jump-diffusion process $(Z_t)$ is given as the strong solution to
the set of equations (\ref{eq:X}), (\ref{eq:Z}).
The following conditions
will be assumed on  the coefficients and the intensity measure.

\begin{assu}
\label{diffcond}
There exists $\theta \in (0, \infty)$ and a measurable function
$\rho: E \to [0, \infty)$ such that
$$\int_E \rho^2(z) F(dz) < \infty, $$
and the following conditions hold.
\begin{description}
\item{(i)} Lipschitz Condition:  For all $y, y' \in \R^n$, $z \in E$,
$$
|\beta(y) - \beta(y')| + |a(y) - a(y')| \le \theta |y-y'|,$$
$$
|h(\delta(y,z)) - h(\delta(y',z))| \le \rho(z) |y-y'|,$$
$$
|h'(\delta(y,z)) - h'(\delta(y',z))| \le \rho^2(z) |y-y'|.$$
\item{(ii)} Growth Condition:  For all $y \in \R^n$, $z \in E$,
$$
\frac{|\beta(y)|}{1 + |y|} + |a(y)| \le \theta, $$
$$
|h(\delta(y,z))| \le \rho(z),$$
$$
|h'(\delta(y,z))| \le \rho^2(z) \wedge \rho^4(z).
$$
\end{description}
\end{assu}

Under the above conditions it can be shown that
there is a unique strong solution to (\ref{eq:X})
and (\ref{eq:Z}) which is a strong Markov process.
I.e., the following result holds.
\bt
\label{exuni}
Suppose that Conditions  \ref{regular} and \ref{diffcond} hold, and
that on $(\Om,\F,\F_t,P)$ we are given processes $(W,N)$ as
above.
Then, for all $x \in G$ there exists, on the basis $(\Om,\F,\F_t,P)$,
 a unique pair of
$\{\F_t\}$-adapted processes
$(Z_t, k_t)_{t \ge 0}$ with paths in $D([0,\infty):\R^n)$, and a
progressively measurable process $(\gamma_t)_{t \ge 0}$, such that
the following hold:
\begin{enumerate}
\item
$Z_t \in G$, for all $t \ge 0$, a.s.
\item  For all $t \ge 0$,
\beqn
Z_t&=&x+\int_0^t\beta(Z_s)ds+\int_0^ta(Z_s)dW_s
  +\int_{[0,t]\times E}h(\del(Z_{s-},z))[N(ds,dz)-q(ds,dz)]\non \\
 &+& \int_{[0,t]\times E}h'(\del(Z_{s-},z))N(ds,dz) +
k_t,
\label{refdiff}\eeqn
a.s.
\item
For all $T \in [0, \infty)$
$$
|k|_T < \infty , \;\; a.s.
$$
\item
$$|k|_t = \int_0^t I_{\{Z_s \in \pa G\}} d|k|_s,$$
and
$k_t = \int_0^t \gamma_s d|k|_s$ with $\gamma_s \in
d(Z_s)$
a.e. $[d|k|]$.
\end{enumerate}
Furthermore, the pair $(Z_t-k_t, Z_t)$
is the unique $\{\F_t\}$-adapted pair of processes
with cadlag paths which satisfies equations (\ref{eq:X}, \ref{eq:Z})
for all $t$, a.s., with the given driving terms $(W,N)$.
Finally, $(Z_t)$ is a strong Markov process on $(\Om,\F,\F_t,P)$.
\end{theorem}
The proof of the theorem follows via the usual Picard iteration method
on using the Lipschitz property of the SM.  We refer the reader to
\cite{dupish1} where a similar argument for constrained diffusion
processes is presented.
\begin{remark}\rm
Condition \ref{diffcond} is a version of the assumptions in
\cite{jacshi}, Chapter III, where strong existence and uniqueness
results for unconstrained jump-diffusion processes are considered.
The conditions assumed there are
substantially weaker, and can be similarly weakened in the current
context as well, via similar arguments.
\end{remark}
\begin{remark}\rm
Taking $a(z) \equiv a$, $\beta(z) \equiv \beta$, $(E, \cle) \equiv
(\R^n, \clb(\R^n))$, $\delta(y,z) \Df z$, $\rho(z) \Df |z|1_{|z| \le
1} + \sqrt{|z|}1_{|z| \ge 1}$ and $F(dz) \equiv K(dz)$, we see that a
L\'evy process satisfying Condition \ref{finmean} is a special case of
the process $\{Z_t\}$ in Theorem \ref{exuni}.
\end{remark}

Here are the main results of this paper.
The first result gives sufficient conditions for
transience and stability of a reflecting L\'evy process.
The transience proof
is a simple consequence of the law of large numbers, while the
stability is treated in a more general framework in the context
of a reflected jump-diffusion process.
For a Borel set $A \subset G$, let $\tau_A$
denote the first time $Z$ hits $A$.
Define
\be
\label{condefa}
 \C \Df \cone\{ - d_i, i\in \{1, \ldots, N\} \}. \ee
\begin{theorem}
\label{th:levy} Let $X$ and $Z$ be as in Definition
\ref{reflevy}.
Assume that Conditions \ref{compS}, \ref{regular} and \ref{finmean} hold.
\ben
\item If $EX_1\in\calC^c$, then there is a constant
$\gamma\in G\setminus\{0\}$
such that for all $x \in G$, $Z_t/t\to \gamma$ as $t\to\infty$,
$P_x$-a.s.
\item
If $EX_1\in\calC^o$, then there  is a
compact set $A$ such that for all $M \in (0, \infty)$,
$$
 \sup_{z \in G, |z| \le M}E_z\tau_A < \infty.
$$
\een
\end{theorem}
Next we consider reflected jump-diffusion processes.
If in equation (\ref{eq:X}) $X$ were replaced by $Z$, and the coefficients
$a$, $\beta$ and $\del$ were extended to all of $\R^n$, then this equation
alone would define a diffusion process with jumps $Z$, the
extended generator of which we denote by $\calL$
(see \cite{jacshi}, Chapter IX, p. 514 for the form of the extended
generator in this setting). Let
$\id:\R^n\to\R^n$ denote the identity map, and define
\be \label{neti}
\til\beta \doteq\calL\,\id=\beta+\int_{E}h'(\del(\cdot,z))F(dz).
\ee
Note that in view of Condition \ref{diffcond}, there is
a constant $\al < \infty$ such that
\be\label{bound:int} \sup_{x \in
\R^n} \lt|\int_{E}h'(\del(x,z))F(dz)\rt|\le \al.
\ee
We use the
generator of the ``unconstrained'' jump-diffusion process only as a
motivation to define the vector field $\til\beta$. Since we only deal with
constrained diffusions, we will consider only the restriction of
$\til\beta$ to $G$, which, with an abuse of notation, we still
denote by $\til\beta$.
Of course, $\tilde\beta$ can otherwise be defined by
the right hand side of (\ref{neti}).
Our main assumption on $\til\beta$ is the
following.
\begin{assu}\label{cond:cone2}
There exists a compact set $\calC_1$ contained in the
interior $\calC^o$ of $\calC$ such that the
range of $\til\beta$ is contained in $\cup_{k \in \N}\, k\calC_1$.
\end{assu}
Here is the main
result on the stability of reflected jump-diffusions.
\begin{theorem}\label{th:rec}
Let $(Z_t)$ be as in Theorem \ref{exuni}.  Suppose that
Conditions \ref{compS}, \ref{regular}, \ref{diffcond} and
\ref{cond:cone2} hold. Then there is a
compact set $A$ such that for any compact $K\subset G$,
$\sup_{z\in K}E_z\tau_A < \infty$.
\end{theorem}
\begin{remark}\rm \label{expbdsep}
We will, in fact, obtain a more precise bound, namely
$E_z\tau_A\le \al|z|+1$, for some constant $\al$ independent of $z\in G$.
\end{remark}
As an immediate corollary of the above theorem we have the following result.
\begin{cor}
\label{sspace}
Let $p(t,x,dy)$ denote the transition probability function of the
Markov process $\{Z_t\}$.  Suppose that there is a closed set $S
\subset G$ such that $p(t,x, S) = 1$ for all $x \in S$ and $t \in (0,
\infty)$.  Let the compact set $A$ be as in Theorem \ref{th:rec} and
suppose that the assumptions of that theorem hold.  Then $\sup_{z \in
S, |z| \le M} E_z(\tau_{A \cap S}) < \infty$ for all $M \in (0,
\infty)$.
\end{cor}
The following result on ``boundedness in probability'' of the process
$\{Z_t\}$ is a consequence of the existence of a suitable Lyapunov
function and will be proved in Section 3.
\bt
\label{bdinp}
Let the assumptions of Theorem \ref{th:rec} hold.
Then for every $M \in (0, \infty)$, the family of
probability measures,
$\{ P_{z}( Z(t) \in \cc);\; t \in [0, \infty),\; z \in G \cap B(0,M)\}$
is tight.
\et
 From the above result we have, on using the Feller property of
$\{Z_t\}$, the following corollary.
\begin{cor}
\label{oneinv}
Suppose the assumptions of Theorem \ref{th:rec} hold.
Then the Markov process
$\{Z_t\}$ admits at least one invariant measure.
\end{cor}
We now impose the following communicability condition
on the Markov process $\{Z_t\}$ relative to a set $S$.
\begin{assu}
\label{commune}
Let $S$ be as in Corollary \ref{sspace} and let $\nu$ be a
$\sigma$-finite measure with support $S$.  Then for all $r \in (0,
\infty)$ and $C \in \clb(\R^n)$ with $\nu(C) > 0$, $\inf_{x \in S,
|x|\le r} P_x(Z_1 \in C) > 0$.
\end{assu}
The above assumption is satisfied with $S = G$ and $\nu$ as the Lebesgue
measure, if the
diffusion coefficient $a$ in (\ref{refdiff}) is uniformly non degenerate.

Now we can give the following result on positive Harris recurrence.
The proof of the theorem is similar to that of Theorem 2.2 of \cite{abd} and
thus is omitted.
\bt
\label{harris}
Let the assumptions of Theorem \ref{th:rec} and Corollary \ref{sspace}
hold.  Further suppose that Condition \ref{commune} holds.  Then for
all closed sets $C$ with $\nu(C) > 0$, and all $M>0$,
we have that $\sup_{z \in S,|z| \le M} \E_z(\tau_C) < \infty .$
\et
Finally, we introduce one more
condition which again is satisfied if the diffusion coefficient is
uniformly non degenerate and $\nu$ is the Lebesgue measure on $G$.
\begin{assu}
\label{resolve}
For some $\lambda \in (0, \infty)$, the probability measure $\theta$ on $G$
defined as
$$
\theta(F) \Df \lambda \int_0^{\infty} e^{-\lambda t} p(t,x,F) dt, \;\;
F \in \clb(G)
$$
is absolutely continuous with respect to
$\nu$.\end{assu} The following theorem is a direct consequence of
Theorem 4.2.23, Chapter 1 of \cite{skor}.

\bt \label{uniqinv}
Let the assumptions in Theorem \ref{harris}
hold.  Further suppose that Condition \ref{resolve} holds.
Then $(Z_t)$ has a unique invariant probability measure.
\et

\section{Proofs of the main results}
\beginsec
We begin with the proof of Theorem \ref{th:levy}.

\noindent{\bf Proof of Theorem \ref{th:levy}:} Since part 2 is a special
case of Theorem \ref{th:rec} (note that Condition \ref{diffcond}
implies, in the special case of a L\'evy process, Condition \ref{finmean}),
we consider only part 1. Let
$\beta=EX_1$ and let $\phi(t)=\beta t$. Then by \cite[Lemma 3.1
and Theorem 3.10(2)]{BD}, $\Gamma(\phi)(t)=\gamma t$, where
$\gamma\ne0$. By the Lipschitz continuity of $\Gamma$,
$Z_t=\Gamma(\phi)(t)+\la_t=\gamma t+\la_t$, where
$$
 |\la_t|\le\ell\sup_{s\le t}| z_0 +X_s-s\beta|.
$$
 From the strong law of large numbers  $t^{-1}(X_t-t\beta)\to 0$ a.s.
Combined with a.s.\ local boundedness of $X$, this implies that
$\sup_{s\le t}|z_0+X_s-s\beta|\to0$ a.s.
Thus $t^{-1}|\la_t|\to0$ a.s., and this proves the result. \ink

In the rest of the paper we prove Theorem \ref{th:rec} and
its consequences. Hence we will assume throughout that Conditions
\ref{compS}, \ref{regular}, \ref{diffcond} and \ref{cond:cone2} hold.
The proof of Theorem \ref{th:rec} is based on the existence of a
suitable Lyapunov function which is defined as follows.
\begin{definition}\label{def:lyap}
\ben
\item
We say that a function $f\in C^2(G\setminus\{0\})$ is a {\em
Lyapunov function for the SP $(G,d)$ with respect to the mean
velocity $r_0$,} if the following conditions hold.
\ben
\item
For all $N \in (0, \infty)$, there exists $M \in (0, \infty)$
 such that  $(x\in G,\ |x|\ge M)$ implies
that $f(x)\ge N$.
\item
For all $\eps>0$ there exists $M \in (0, \infty)$ such that $(x\in G,\
 |x|\ge M)$ implies $\| D^2f(x)\|\le\eps$.
\item
There exists $c \in (0, \infty)$ such that $Df(x)\cdot r_0 \le -c$,
$x\in G\setminus\{0\}$, and $Df(x)\cdot d\le-c$, $d\in d(x)$, $x\in\pl
G\setminus\{0\}$.
\item
There exists $L \in (0, \infty)$ such that
$\sup_{x \in G} |Df(x)| \le L$.
\een
\item
We say that a function $f\in C^2(G\setminus\{0\})$ is a {\em
Lyapunov function for the SP $(G,d)$ with respect to the set
(of mean velocities)
$\til R\subset\R^n$,} if it is a Lyapunov function for the SP
$(G,d)$ with respect to the mean velocity $r_0$, for any
$r_0\in\til R$, and if in item (c) above, the constant $c$ does
not depend on $r_0\in \til R$. \een
\end{definition}
\begin{remark}\rm \label{morem}
{\bf (a)}
If $f$ is a Lyapunov function for the SP $(G,r)$ with
respect to a certain set $\til R$, then $Df$ is Lipschitz
continuous on $\{x\in G:|x|\ge M\}$  with parameter $\eps$, where
$\eps>0$ can be taken arbitrarily small by letting $M$ be large.
This implies a useful consequence of
 the second part of item
(c) in Definition \ref{def:lyap} as follows: There exist
$M_0 \in (0, \infty)$, $\delta_0 \in (0,1)$ such that
$Df(x)\cdot d\le -{c}/{2}$,
whenever $d\in d(y),\; |y-x|\le\del_0,\; y\in \partial G,\; |x|\ge M_0$.

\noi{\bf(b)} If $f$ is a Lyapunov function for $(G,d)$ with
respect to a set $\til R$, then it is automatically a Lyapunov
function for $(G,d)$ with respect to $\cup_{k\in\N}\,k\til R$.
\end{remark}
We say that a function $f$ is radially linear
on $G$ if $f(s x) = s f(x)$ for all $s \in (0, \infty)$
and $x \in G$.

The following result is key to the proof of Theorem \ref{th:rec}
and will be proved in Section 4.
\bt
\label{lyexsts}
Let the assumptions of Theorem \ref{th:rec} hold.
Then there exists a Lyapunov function $f$ for the SP $(G, d)$
with respect to the set $\calC_1$, where $\calC_1$ is as in
Condition \ref{cond:cone2}.  Furthermore, $f$ is radially linear
on $G$.
\et

We now turn to the proof of Theorem \ref{th:rec}.
Write (\ref{refdiff}) as
$$
 Z_t=z_0+\int_0^t\til\bet(Z_s)ds + \int_0^t a(Z_s)dW_s
 +M_t^{(1)}+M_t^{(2)} + k_t,
$$
where
$$
 M_t^{(1)}=\int_{[0,t]\times E}h(\del(Z_{s-},z))[N(ds,dz)-q(ds,dz)],
$$
$$
 M_t^{(2)}=\int_{[0,t]\times E}h'(\del(Z_{s-},z))[N(ds,dz)-q(ds,dz)]
$$
and
$
\ti \beta(\cc)$ is as in (\ref{neti}).
Note that the term that has been subtracted and added is finite
(e.g., by (\ref{bound:int})). Let also
$$
 U_t=\int_0^t\til\beta(Z_s)ds,
$$
and
$$
 M_t=\int_0^t a(Z_s)dW_s+M_t^{(1)}+M_t^{(2)}.
$$
Then \be \label{eq:ZU}
 Z_t=z_0+U_t+M_t+k_t.
\ee
Let $f$ be as in Theorem \ref{lyexsts}.
From Condition \ref{cond:cone2}, it follows (see also Remark
\ref{morem}(b)) that
 \be \label{eq:U} Df(x)\cdot \dot U_t\le -c, \quad x\in
G\setminus\{0\}, \quad t\ge0, \ee
where $c$ is as in Definition \ref{def:lyap}.
For any $\kap \in (0, \infty)$
and compact set $A \subset \R^n$,
 define the sequences $(\tilde\sig_n)$, $(\sig_n)$ of
stopping times as $\til\sig_0=0$,
$$
 \til\sig_n=\til\sig_n(\kap)
 \Df\inf\{t>\til \sig_{n-1}:|X_t-X_{\ti\sig_{n-1}}|\ge \kap\},
$$
$$
 \sig_n=\sig_n(\kap,A)\Df \til\sig_n\w \tau_A.
$$
Let also
$$
 \til n(t)=\inf\{n:\til\sig_n\ge t\},
$$
$$
 \bar n(t)=\bar n(t,A)=\inf\{n:\sig_n\ge \tau_A\w t\},
$$
where the infimum over an empty set is $\infty$.
Note that $\til n(t\w \tau_A) =\bar n(t)$, a.s.
The following are the main lemmas used in the proof of
Theorem \ref{th:rec}.
\begin{lemma}\label{lem:Mt}
$\{M^{(1)}_t\}$ and $\{\int_0^t a(Z_s) dW_s\}$ are
square integrable martingales and $\{M^{(2)}_t\}$ is a  martingale.
\end{lemma}

\begin{lemma}\label{lem:Nt}
There exists a constant $c_1=c_1(\kap) \in (1, \infty)$ such that for
any bounded stopping time $\tau$, $E\til n(\tau)\le c_1 (E\tau + 1)$.
\end{lemma}
For $s \in [0, \infty)$, and a cadlag
process $\{Y_t\}$, we write $Y_s -Y_{s-}$ as $\Delta Y_s$.
\begin{lemma}\label{lem:DelX}
There is a $b_0 \in (0, \infty)$ and a function $\bar\al: [0, \infty) \mapsto [0, \infty)$
with $\bar\al(b)\to0$ as $b\to\infty$
such that for any bounded stopping time $\tau$, and $b > b_0$,
$$
 E\sum_{s\le\tau}|\Del X_s|1_{|\Del X_s|>b} \le\bar\al(b)E\tau.
$$
\end{lemma}

\noindent{\bf Proof of Theorem \ref{th:rec}:}
Let $f$ be as in Theorem \ref{lyexsts} and let $a_0$
be as in (\ref{asn:1}). For any $M$, the level set $\{x\in
G:f(x)\le M\}$ is compact. If $\til M=\til M(M)=\max\{|x|:f(x)\le
M\}$, then the set
\be \label{setA}
 A=A(M)=\{x\in G:x\cdot a_0\le \til M\}
\ee contains the level set, and is compact. In addition,
$G\setminus A$ is convex. Definition \ref{def:lyap}.1(a),(b)
implies that there is a function $\eps_0(\til M)$ such that $\|
D^2f(x)\|\le\eps_0(\til M)$ for $x\in G\setminus A(M)$, and where
$\til M=\til M(M)$ is as above, and $\eps_0(\til M)\to 0$ as
$\til M\to\infty$. The notation $\eps_0(\til M)$ and $\til M(M)$
is used in what follows.

Write
$x_n=X_{\sig_n}$, $x_{n-}=X_{\sig_n-}$, where $\{X_t\}$
is as in (\ref{eq:Z}). Define
similarly $k_n$, $k_{n-}$, $z_n$, $z_{n-}$, $u_n$, $u_{n-}$,
$m_n$, $m_{n-}$ for the processes $k,Z,U$ and $M$, respectively. Let
$\kap$ be so small that $2\ell\kap\le\del_0/2$, where $\del_0$ is
as in Remark \ref{morem}(a). $\kap$ will be fixed throughout. The
proof will be based on establishing a bound on
$Ef(z_m)=f(z_0)+\sum_1^mE[f(z_n)-f(z_{n-1})]$. According to
(\ref{eq:ZU}), one has
\beaa
 z_n-z_{n-1} &=& x_n-x_{n-1}+k_n-k_{n-1} \\
 &=&  u_n-u_{n-1}+m_n-m_{n-1}+k_n-k_{n-1}.
\eeaa
We consider two cases.

\noindent{\em Case 1:} $|x_n-x_{n-1}|\le2\kap$.

Consider the linear interpolation $z^\theta$ defined for
$\theta\in[0,1]$ as
$$
 z^\theta = z_{n-1}+\theta(z_n-z_{n-1}).
$$
Then \be\label{eq:delf}
 f(z_n)-f(z_{n-1})=\int_0^1Df(z^\theta)d\theta\cdot(z_n-z_{n-1}).
\ee
By the Lipschitz continuity of the SM,
$z^\theta\in B_{2\ell\kap}(z^0) \subset B_{\delta_0/2}(z^0)$
 for $\theta\in[0,1]$. Also note that
for $s \in [\sigma_{n-1}, \sigma_n]$,
\be\label{eq:Ry}
 \gamma_s \in\bigcup_{x\in B_{2\ell\kap}(z^0)}d(x)
 \subset \bigcup_{x\in B_{\del_0/2}(z^0)}d(x), \quad [d|k|]\; a.s.
\ee
Let $M$ be so large that $\til M\ge M_0+1$,
where $M_0$ is as in Remark \ref{morem}(a). Then any $x\in A=A(M)$ satisfies
$|x|\ge M_0+1$. By convexity of $A$ we therefore have for $n\le
\bar n(\infty)$ that $|z^\theta|\ge\til M\ge M_0+1$,
and we get from (\ref{eq:Ry}) that for $\theta \in [0,1]$,
\beqn
 Df(z^\theta)\cdot[k_n-k_{n-1}] &=& Df(z^\theta) \int_{\sigma_{n-1}}^{\sigma_n}
 \gamma_s d|k|_s \non \\
&\le & -\frac{c}{2}(|k|_n-|k|_{n-1})\le 0.
  \label{eq:DfY}
\eeqn
Now let $\ep \Df \ep_0(\ti M - \ell b)$.
By (\ref{eq:U}),
$Df(z^0)\cdot(u_n-u_{n-1})\le-c(\sig_n-\sig_{n-1})$. Therefore,
from part (b) of Definition \ref{def:lyap}, and (\ref{eq:delf}),
(\ref{eq:DfY}), we have
\bea \label{bound:1} f(z_n)-f(z_{n-1})
&\le&
 Df(z_{n-1})(z_n-z_{n-1})+\int_0^1|Df(z^\theta)-Df(z^0)|d\theta
 \,|z_n-z_{n-1}| \nonumber \\
 &\le&
 Df(z_{n-1})(z_n-z_{n-1}) + (\eps)(2\ell\kap)(2\ell\kap) \nonumber \\
 &\le&
 -c(\sig_n-\sig_{n-1})+Df(z_{n-1})(m_n-m_{n-1}) + 4(\ell)^2\eps\kap^2.
\eea

\noindent{\em Case 2:} $|x_n-x_{n-1}|>2\kap$.

The argument applied in Case 1 gives an analogue of
(\ref{bound:1}) in the form \be\label{bound:2}
f(z_{n-})-f(z_{n-1}) \le
 -c(\sig_n-\sig_{n-1})+Df(z_{n-1})(m_{n-}-m_{n-1}) + 4\ell\eps\kap^2.
\ee
Next we provide a bound on $f(z_n)-f(z_{n-})$.
Let
$$
 \hat z^\theta=z_{n-}+\theta(z_n-z_{n-}), \quad \theta\in[0,1].
$$
Note that $k_n-k_{n-}\in d(z_n)$, by Definition \ref{def-sp},
and therefore $Df(z_{n})\cdot(k_n-k_{n-})\le0$. Also,
recall that $|Df|\le L$. Let $b_0$ be as in Lemma \ref{lem:DelX} and $b> b_0$ be arbitrary. Then if
$|x_n-x_{n-}|\le b$, then for all $\theta\in[0,1]$, $|\hat
z^\theta-\hat z^0|\le \ell b$, and therefore the bound $\|D^2
f(z^\theta)\|\le\eps\doteq\eps_0(\til M(M)-\ell b)$
holds.
Thus
\bea\label{bound:3}
f(z_n)-f(z_{n-}) &=& \nonumber
 \lt(\int_0^1 Df(\hat z^\theta)d\theta\rt)\cdot[(x_n-x_{n-})+(k_n-k_{n-})]\\
 &=& \nonumber
 Df(z_{n-})\cdot (x_n-x_{n-})+\int_0^1(Df(\hat z^\theta)-Df(z_{n-}))d\theta
 \cdot (x_n-x_{n-}) \\
 && \nonumber
 +Df(z_n)\cdot (k_n-k_{n-})+\int_0^1(Df(\hat z^\theta)-Df(z_n))d\theta
 \cdot (k_n-k_{n-})\\
 &\le& \nonumber
 Df(z_{n-})\cdot(x_n-x_{n-})+2\eps\ell^2b^2
 +4\ell L|x_n-x_{n-}|1_{|x_n-x_{n-}|>b}\\
 &=&
 Df(z_{n-})\cdot(m_n-m_{n-})+2\eps\ell^2b^2
 +4\ell L|x_n-x_{n-}|1_{|x_n-x_{n-}|>b}.
\eea

Let $J_q(t)$ denote the set $\{n\le\bar
n(t)\w q:|x_n-x_{n-1}|>2\kap\}$. Combining (\ref{bound:1}),
(\ref{bound:2}) and (\ref{bound:3}) we get \beaa \sum_{n\le\bar
n(t)\w q}f(z_n)-f(z_{n-1}) &\le& -c\sig_{\ov n(t) \w q}+
\sum_{n\le\bar n(t) \w q} Df(z_{n-1})(m_n-m_{n-1})
 +4\ell\eps\kap^2\ov n(t) \w q \\
 &+&
 \sum_{n\in J_q(t)}(Df(z_{n-})-Df(z_{n-1}))\cdot(m_n-m_{n-}) \\
 &+&
 \sum_{n\in J_q(t)}(2\eps\ell^2 b^2+4\ell L|x_n-x_{n-}|1_{|x_n-x_{n-}|>b}).
\eeaa Since $m_n-m_{n-}=x_n-x_{n-}$, the following inequality
holds
$$
 \sum_{J_q(t)}(Df(z_{n-})-Df(z_{n-1}))\cdot(m_n-m_{n-})
 \le \sum_{J_q(t)}(2\eps\ell^2 b^2+4\ell L|x_n-x_{n-}|1_{|x_n-x_{n-}|>b}).
$$
Writing $\til M_t=\sum_{n\le\bar n(t)\w q}Df(z_{n-1})(m_n-m_{n-1})$,
we get
\beqn
 f(z_{\ov n(t) \w q})-f(z_0) &=&
 \sum_{n\le \bar n(t)\w q}(f(z_n)-f(z_{n-1})) \non \\
 &\le&
 -c(\sigma_{\ov n(t) \w q})+\til M_t +4\ell\eps\kap^2\bar n(t)
 \non \\
 &+&
 2\eps\ell^2 b^2|J_q(t)|
 +\sum_{n\in J_q(t)}4\ell L|x_n-x_{n-}|1_{|x_n-x_{n-}|>b}. \label{sep1}
\eeqn
From Lemma \ref{lem:Mt},  $E(\til M_t) = 0$.
Using Lemma \ref{lem:Nt}
$$E(\bar n(t)) = E(\ti n(t \w \tau_A) \le c_1 (E(t \w\tau_A)+1).$$
Observing that
$$
\sum_{n\in J_q(t)}|x_n-x_{n-}|1_{|x_n-x_{n-}|>b} \le
\sum_{s \le t \w \tau_A} |\Delta X_s| 1_{|\Delta X_s| > b},$$
we have from Lemma \ref{lem:DelX} that
the expectation of the term on the left side above
is bounded by $\bar\alpha(b) E(\tau_A \w t)$.
Combining these observations, we have that
\beaa
 Ef(z_{\ov n(t) \w q}) -f(z_0) &\le&
 -cE(\sigma_{\ov n(t) \w q})
 +\eps(4\ell\kap^2+2\ell^2 b^2)(c_1(E(\tau_A\w t) + 1)\\
 &+&
 4\ell L\bar\al(b)E(\tau_A\w t).
\eeaa Let $b$ be so large that $4\ell L\bar\al(b)\le c/3$.
Recalling the definition of $\eps$,
let $M$ be so large, thus $\eps$
so small, that $\eps(4\ell\kap^2+ 2\ell^2 b^2)c_1\le c/3$. Then
$$
-f(z_0) - c/3\le -cE(\sigma_{\ov n(t) \w q}) +\frac {2c} 3 E(\tau_A\w t),
$$
Taking $q \to \infty$, and recalling that $\sigma_{\ov n(t)} \ge \tau_A \w t$,
we see that
 $E(\tau_A\w t)\le 3f(z_0)/c + 1$. Finally,  taking $t\to\infty$, we
get for each $z_0$, $E_{z_0}\tau_A\le3f(z_0)/c+1$.
Note that $\kap,c_1,c,L$ do not
depend on $z_0$, nor do the choices of $M,\eps(M),b,\bar\al(b)$. The
result follows. \ink

We now present the proof of Theorem \ref{bdinp}.\\
{\bf Proof of Theorem \ref{bdinp}.}
The proof is adapted from \cite{kushweak}, pages 146-147.
Since the Lyapunov function $f$ satisfies
$f(z) \to \infty$ as $|z| \to \infty$, it suffices to show that:
For all $\delta > 0$ and $L_0 \in (0, \infty)$
there exists an $\eta$ such that
\be
\label{tigffn}
\inf_{x \in G, |x| \le L_0} P_x(f(Z(t)) \le\eta) \ge 1- \delta .
\ee
Let $A$ be as in the proof of Theorem \ref{th:rec}.  Fix $\la > \ti M$
and define
$A_{\la} \Df \{x \in G: x\cdot a_0 \ge \la\}$.
Let $\ov \la \Df \sup\{|x|:x \in A_{\la}^c\}$ and set
$\rho \Df \sup\{f(x): x \in G \cap B(0,1)\}$.
Recalling the radial property
of the Lyapunov function we have that for all $x \neq 0$;
$f(x) \le \rho |x|$.

Now we define a sequence of stopping times $\{\tau_n\}$ as follows.
Set $\tau_0 = 0$. Define
$$
\tau_{2n+1} \Df \inf \{t > \tau_{2n}: Z(t) \in A\}; \;\; n \in \N_0$$
and
$$
\tau_{2n+2} \Df \inf \{t > \tau_{2n+1}: Z(t) \in A_{\la}\}; \;\; n \in \N_0.$$
Without loss of generality, we assume that $\tau_n < \infty$ with
probability $1$ for all $n$.
From Remark \ref{expbdsep} we have that, for all $n \in \N_0$,
\beqn
E(\tau_{2n+1} - \tau_{2n} \mid \clf_{\tau_{2n}})
& \le & c |Z_{\tau_{2n}}| + 1 \non \\
& \le & c(|\Del Z_{\tau_{2n}}|+\bar\la) + 1 \non \\
& \le & \alpha |\Delta Z_{\tau_{2n}}| + \alpha. \label{sep3a}
\eeqn
Next observe that, for all $\eta > \ov \la \rho$:
\be \label{sep3b}
f(z(t))\le\eta, \quad t\in[\tau_{2n+1},\tau_{2n+2}), n\in\N_0.
\ee
Now we claim that there is a constant $\al$
such that for all $\eta > 0$ and $x \in G$
\be \label{sepfin2}
P_x(\sup_{0 \le t < \tau_{1}} f(Z(t)) \ge \eta)
 \le \al \frac{f(x) + 1}{\eta}\ee
and for all $n \in \N$
\be \label{sepfin1}
P(\sup_{\tau_{2n} \le t < \tau_{2n+1}} f(Z(t)) \ge \eta
\mid \clf_{\tau_{2n}}) \le \al \frac{|\Delta Z_{\tau_{2n}}| + 1}{\eta}.\ee
We only show (\ref{sepfin1}), since the proof of (\ref{sepfin2}) is similar.
By arguing as in the proof of Theorem
\ref{th:rec} (see (\ref{sep1})), we have that
\be
\label{sep2a}
\sup_{\tau_{2n} \le t < \tau_{2n+1}} f(Z(t))
\le f(Z_{\tau_{2n}}) + L\ell \kappa
+ \sup_{1 \le k \le \ov n(\tau_{2n+1})} \sum_{1 \le j \le k}
(f(z_{j}) - f(z_{j-1}))
\ee
where $\{\sigma_j\}$ and $\ov n(\cc)$ are defined as in the displays
below (\ref{eq:U}) with $\ti \sigma_0 \Df \tau_{2n}$ (rather than $0$) and
$\tau_A$ replaced by $\tau_{2n+1}$.
Given a stopping time $\;\tau$, denote the
conditional expectation and conditional probability with respect to the $\sigma$-field
$\clf_{\tau}$ by $\E_{\tau}$ and $\IP_{\tau}$ respectively.  Then, we have  via arguments as in
Theorem \ref{th:rec} that
\begin{eqnarray*}
&&
\E_{\tau_{2n}}(\sup_{1 \le k \le \ov n(\tau_{2n+1})} \sum_{1 \le j \le k}
(f(z_{j}) - f(z_{j-1})) )\\
&&
\le \E_{\tau_{2n}}(\sup_{1 \le k \le \ov n(\tau_{2n+1})}
|\sum_{1 \le j \le k} Df(z_{j-1})(m_j - m_{j-1})|)\\
&&
\qquad + \al( \E_{\tau_{2n}}(\tau_{2n+1} - \tau_{2n}) + 1).
\end{eqnarray*}
Doob's inequality yields that
$$
\E_{\tau_{2n}}(\sup_{1 \le k \le \ov n(\tau_{2n+1})}
|\sum_{1 \le j \le k} Df(z_{j-1})(m_j - m_{j-1})|)
\le \al(\E_{\tau_{2n}}(\tau_{2n+1} - \tau_{2n}) + 1).
$$
Combining the above observations with (\ref{sep3a}) we have that
$$
\E_{\tau_{2n}}(\sup_{1 \le n \le \ov n(\tau_{2n+1})} \sum_{1 \le j \le n}
(f(z_{j}) - f(z_{j-1})) ) \le \al (|\Delta Z_{\tau_{2n}}| + 1).
$$
Combining this with (\ref{sep2a}) we have (\ref{sepfin1}).

Following \cite{kushweak} we can choose an
integer $k_{\delta}$ and, for each $t$, an integer valued random variable
$j(t, \delta)$ such that
$\tau_{j(t,\delta)}$ are stopping times and
$$P(\tau_{j(t,\delta)} \le t \le
\tau_{j(t,\delta)+ k_{\delta}}) \ge 1 - \delta/2.$$
Now define $J_i \Df [\tau_{j(t,\delta)+ i-1}, \tau_{j(t,\delta)+ i})$
and fix $\eta > \ov \la \rho$.
Let  $\tau'$ be the hitting time of the set $A_{\la}$ by
$Z_t$.
Then
\beqn
P_x(f(Z(t)) \ge \eta)
& \le & \frac{\delta}{2} +
\sum_{i=1}^{k_{\delta}} P_x(\sup_{s \in J_i} f(Z(s)) \ge \eta) \non \\
& \le &
\frac{\delta}{2}
+ \sum_{i=1}^{k_{\delta}} E_x(\E_{\tau_{j(t,\delta)}}
(\IP_{\tau_{j(t,\delta)+ i-1}}(\sup_{s \in J_i} f(Z(s)) \ge \eta))) \non \\
&+& P_x( \sup_{0 \le s \le \tau_1} f(X(s)) \ge \eta)\non \\
& \le &
\frac{\delta}{2} + \frac{\al}{\eta}
\sum_{i=1}^{k_{\delta}} E_x(\E_{\tau_{j(t,\delta)}}
(|\Delta Z_{\tau_{j(t,\delta) + i-1}}| + 1)) \non \\
&+&
\frac{ \al(f(x) +1)}{\eta}\non \\
&\le&
\frac{\delta}{2} + \frac{\al(f(x) +1 + k_{\delta} + b k_{\delta})}{\eta}\non \\
&+& \frac{\al}{\eta} E_x( \E_{\tau_{j(t,\delta)}}
(\sum_{i=1}^{k_{\delta}}|\Delta Z_{\tau_{j(t,\delta) + i-1}}|1_{
|\Delta Z_{\tau_{j(t,\delta) + i-1}}| > b})),  \label{ramipf}
\eeqn
where the third inequality above is the consequence of (\ref{sep3b}),
(\ref{sepfin1}) and  (\ref{sepfin2}) and in the fourth inequality $b \in (1, \infty)$
is arbitrary.
Next note that
\beq
&&E_x( \E_{\tau_{j(t,\delta)}}
(\sum_{i=1}^{k_{\delta}}|\Delta Z_{\tau_{j(t,\delta) + i-1}}|1_{
|\Delta Z_{\tau_{j(t,\delta) + i-1}}| > b}))\\
&\le &
E_x( \E_{\tau_{j(t,\delta)}}
(\sum_{s \in [\tau_{j(t,\delta)}, \tau_{j(t,\delta)+k_{\delta} + 1})}
|\Delta Z_{s}|1_{
|\Delta Z_{s}| > b}))\\
&\le&
E_x( \E_{\tau_{j(t,\delta)}}
(\int_{[\tau_{j(t,\delta)}, \tau_{j(t,\delta)+k_{\delta} + 1}) \times E}
h'(\delta(Z_{s-},z)) 1_{|h'(\delta(Z_{s-},z))| > b} F(dz) ds )) \\
&\le &
(k_{\delta} +1) \bar\alpha(b),
\eeq
where $\bar\alpha(b) \Df \int_{E} \rho^2(z) 1_{\rho^2(z) > b} F(dz)$.

Using the above observation in (\ref{ramipf}) we have that
$$
P_x(f(Z(t)) \ge \eta) \le
\frac{\delta}{2} + \frac{\al(f(x) +1 + k_{\delta} + b k_{\delta})}{\eta}
+ \frac{\al}{\eta} (k_{\delta} +1) \bar\alpha(b).
$$
The result now follows on taking $\eta$ suitably large. \ink

We now give the proofs of the lemmas.

\noi{\bf Proof of Lemma \ref{lem:Mt}:}
Since $a(\cc)$ is a bounded function we have that $\int_0^t a(Z_s) dW_s$
is a square integrable martingale.
In order to show that $M^{(2)}_t$ is a martingale, it suffices
to show, in view of Theorem II.1.8 of \cite{jacshi}
that for all $T \in [0, \infty)$,
$$
\int_{[0,T]\times E} E|h'(\delta(Z_{s-}, z))| q(ds,dz) < \infty .$$
(The cited theorem states a local martingale property, however the proof there
shows the above stronger assertion.)
The  inequality follows on observing that from Condition \ref{diffcond}
the above expression is bounded  by
$
T \int_{E} \rho^2(z) F(dz)<\infty.$
Finally, in view of Theorem II.1.33 of \cite{jacshi}, to show that
$M^{(1)}_t$ is a square integrable martingale, it suffices to show that
$$
\sup_{s \in [0,T]} \int_{E} E|h(\delta(Z_{s-}, z))|^2 F(dz) < \infty .$$
The last inequality follows, once more from Condition \ref{diffcond}.
\ink

\noi{\bf Proof of Lemma \ref{lem:Nt}:}
Recall that
$X_t=z_0+\int_0^t\til\beta(Z_s)ds+\int_0^ta(Z_s)dW_s+
M^{(1)}_t+M_t^{(2)}$, where $M_t^{(i)}$ are martingales.
Since $M^{(1)}_t$ is a square integrable martingale, by
Doob's inequality we have
\beaa
 E\sup_{s\le\eps}|M_s^{(1)}|^2 &\le&
 4E|M_\eps^{(1)}|^2 \\
 &=&
 4E\int_{[0,\eps]\times E}|h(\del(Z_{s-},z))|^2F(dz)ds \\
 &\le &
 4\eps\int_{E}\rho^2(z) F(dz).
\eeaa
Also observe that
\beq
E|M^{(2)}_{\ep}| & \le &
E | \int_{[0,\ep\times E} h'(\delta(Z_{s-}, z)) N(ds, dz) |
+ E| \int_{[0,\ep]\times E} h'(\delta(Z_{s-}, z)) q(ds, dz) |\\
&\le & 2 E \int_{[0,\ep]\times E} |h'(\delta(Z_{s-}, z))| q(ds, dz)\\
& \le & 2\ep \int_{E} \rho^2(z) F(dz),
\eeq
where the second inequality is a consequence of
Theorem II.1.8 of \cite{jacshi} and the last inequality follows from Condition
\ref{diffcond}.
Using the linear growth of $\til\beta$ and the Lipschitz property of
$\Gamma$, the above moment bounds
show that
$$E\sup_{s\le \eps}
|X_s-z_0|\le \al\sqrt{\eps}+\al\int_0^\eps E\sup_{s\le\theta}
|X_s-z_0|d\theta.$$
Hence, by Gronwall's inequality, for every $\del>0$ there
is $\eps>0$ such that $E(\sup_{0 \le s \le\eps} |X_s-z_0|) \le\del$.
By choosing $\eps \in (0,1)$ small enough one can obtain
$$
 P(\til\sig_1\le\eps)=
 P(\sup_{s\le\eps}|X_s-X_0|\ge\kap) \le1/2.
$$

Let $\calF^n=\calF_{\til\sig_n}$. By the strong Markov property
of $Z$ on $\calF_t$, and by considering the martingales
$M_{\til\sig_{n-1}+\eps}^{(i)}-M_{\til\sig_{n-1}}^{(i)}$
in place of $M_\eps^{(i)}$, one obtains that for any $n$,
$P(\til\sig_n-\til\sig_{n-1}>\eps|\calF^{n-1})>1/2$.
    Let $\tau$
be a bounded $(\calF_t)$-stopping time.
An application of Chebychev's inequality
and the observation that since $\ep \in (0,1)$
the sets $\{(\ti \sig_i - \ti \sig_{i-1}) > \ep \}$
and $\{(\ti \sig_i - \ti \sig_{i-1})\w 1> \ep \}$ are equal, we have that
\be
\label{lastone}
 \frac\eps 2 E[\til n(\tau)\w k]\le
 E\sum_{i=1}^{\til n(\tau)\w k}E[(\til\sig_i-\til\sig_{i-1})\w 1|\calF^{i-1}].
\ee
Define
$$
S_j \Df \sum_{i=1}^j \left(
(\til\sig_i-\til\sig_{i-1})\w 1 - E[(\til\sig_i-\til\sig_{i-1})\w 1|\calF^{i-1}]\right)
$$
Then $(S_j, \clf^j)$ is a zero mean martingale.
Observing that
$\til n_\tau$ is a stopping time on the filtration $(\calF^n)$
we have that for all $k \in \N$, $E(S_{\ti n(\tau)\w k}) = 0$.
Hence from (\ref{lastone}) it follows that
$$ \frac\eps 2 E[\til n(\tau)\w k]
\le E\sum_{i=1}^{\til n(\tau)\w k}[(\til\sig_i-\til\sig_{i-1})\w1]
 \le  E(\til\sig_{\ti n(\tau) -1}) + 1  \le E\tau + 1.
$$
Taking $k\uparrow\infty$, the result follows. \ink

\noindent{\bf Proof of Lemma \ref{lem:DelX}:}
Let $b \in (0, \infty)$ be large enough so that
$h(x) = 0$ for $|x| \ge \frac{b}{2}$ and $\sup_{x \in \R^n} |h(x)| \le \frac{b}{2}$.
Now let $\psi: \R^n \to [0, \infty)$ be defined as $\psi(z) \Df |z|1_{b\le|z|\le b'}$,
where $b' \in (b, \infty)$.  Clearly, for all $x \in \R^n$, $\psi(x) = \psi(h'(x))$.
 Now from Theorem II.1.8 of \cite{jacshi} \beaa E \sum_{s\le
\tau}|\Del X_s| 1_{b\le|\Del X_s|\le b'} &=& E\int_{[0,\tau]\times
E}
|h'(\del(Z_{s-},z))| 1_{b\le|h'(\del(Z_{s-},z))|\le b'} F(dz)ds\\
&\le&
E(\tau)\int_{E}\rho^2(z) 1_{b \le \rho^2(z)} F(dz) \\
&\le&  \bar\al(b)E(\tau),
\eeaa
where $\bar\al(b)\to0$ as $b\to\infty$. The result now follows
upon taking $b'\to\infty$. \ink

\section{Construction of the Lyapunov function}

This section is devoted to the proof of Theorem \ref{lyexsts}.
We begin with a stability result on constrained deterministic
trajectories which was proved in \cite{abd}.

Let $\calC_1$ be as in Condition \ref{cond:cone2}.  Let $\del>0$
be such that $\dist(x, \partial \calC) \ge
 \delta$ for all $x \in \calC_1$. Define
$$
 V=\{v\in B:\int_0^t|v(s)|ds<\infty, v(t)\in\calC_1,\, t\in(0,\infty)
 \},
$$
where $B$ is the set of measurable maps $[0,\infty)\to\R^n$.
For $x\in G$ let
$$
 {\bf Z}_x=\{\Gamma(x+\int_0^\cdot v(s)ds):v\in V\}.
$$
\begin{prop}{\bf \cite{abd}}
\label{staba}
\, For any $x\in G$ and $z\in {\bf Z}_x$, the following holds:
$$
|z(t)| \le \frac{\ell^2|x|^2}{\ell|x| + \delta t}, \quad t \in
[0, \infty),
$$
where $\ell$ is the finite constant in (\ref{lip}).
\end{prop}
Using the above result, the following was used in
\cite{abd} as a Lyapunov function:
\be
\label{hittime}
T(x) \Df \sup_z \inf \{t \in [0, \infty): z(t) =
0\},\ee
where the supremum is taken over all trajectories $z\in {\bf Z}_x$.
This function played a key
role in the  proof of positive recurrence of  certain
constrained diffusion processes studied in \cite{abd}.
 The proof in \cite{abd}
uses crucially certain estimates on the exponential moments
of the Markov process.  Since, in the setting of the current
work the Markov process need not even have finite second moment,
the techniques of \cite{abd} do not apply.  However,
we will show that by using the ideas from
\cite{dupwil} and by  suitable smoothing and modifying
 the hitting time function $T(\cc)$, one can obtain  a
Lyapunov function in the sense of Definition \ref{def:lyap}(2)
with $\ti R$ there replaced by $\clc_1$.
Since for $z\in {\bf Z}_x$, $z(s)=0$ implies $z(t)=0$ for $t>s$,
the function $T(\cc)$ can be rewritten as
$$
T(x) \Df \sup_{{\bf Z}_x} \int_0^{\infty} 1_{(0, \infty)}(|z(s)|) ds.
$$
Our first step in the construction is to replace the above
indicator function by a smooth function $\eta$ defined as follows.
Let $\eta: \R \to [0,1]$ be in $\C^{\infty}(\R)$.  Further assume
that
$\eta(z) = 0$ for all $z \in (-\infty, 1]$, $\eta(z) = 1$ for all $z \in
[2,\infty)$
and $\eta'(z) \ge 0$ for all $z \in \R$.
The next step in constructing a $C^2$ Lyapunov function is an
appropriate modification of this new $T(\cc)$ function near the boundary and
a suitable extension of the function to a neighborhood of $G$.

For each $\lambda \in \Lambda$, $\la\ne\emptyset$, fix a vector
$d^{\lambda}$ as in
Condition \ref{compS}.  Define for $\beta , x \in \R^n$
\beq
v(\beta,x) &=& \beta \;\;\; \mbox{for} \; x \in G \\
&=& d^{\lambda(x)}\;\;\; \mbox{for} \; x \not\in G,
\eeq
where $\lambda(x) = \{i \in \{1, \ldots, N\}: \inn{x}{n_i} \le 0 \}.$
Now $\rho \in C^{\infty}(\R^n)$ be such that the
support of $\rho$ is contained in $\{x: |x| \le 1\}$
and $\int_{\R^n}\rho(x) dx = 1$.
Define for $a > 0$
$$v^{a}(\beta , x) \Df
\frac{1}{(a|x|)^n}\int_{\R^n}
\rho(\frac{x-y}{a|x|})v(\beta , y) dy, \;\;\; x \neq 0.$$
Now let $g: \R \to [0,1]$ be a smooth function such that $g(z) =
1$ for $z \in [0, \half]$ and $g(z) = 0$ for $z \in [1, \infty)$.
Define for $i = 1, \ldots, N$, $x \neq 0$ and $\beta \in \R^n$
$$
v_i^a(\beta, x) =
g(\frac{\dist(x, F_i)}{a|x|}) d_i +
\left[1 - g(\frac{\dist(x, F_i)}{a|x|})\right] v^a(\beta, x)$$
and
$$
v_0^a(\beta,x) = g(\frac{\dist(x, G)}{a|x|}) \beta +
\left[1 - g(\frac{\dist(x, G)}{a|x|})\right] v^a(\beta, x).$$

Also set $v_i^a(\beta, 0) = v_0^a(\beta, 0) = {\bf 0}$,
where ${\bf 0} \Df (0, \ldots, 0)'_{1 \times N}$.
Let
$$K^{a}(\beta,x) \Df \conv \{v_i^a(\beta, x); i = 0, 1, \ldots,
N\}.
$$
Finally, define
$$
K^{a}(x) \Df  \cup_{\beta \in \C_1}K^a(\beta,x),\;\;  x \in \R^n . $$

Now we can define our second modification to the hitting time
function.  In this modified form the supremum in (\ref{hittime})
is taken, instead, over all solutions to
the differential inclusion $\dot{\phi}(t) \in K^a(\phi(t)); \;
\phi(0) = x$.  More precisely,
for a given $x \in \R^n$ let $\phi(\cc)$ be an absolutely
continuous function on $[0, \infty)$  such that
$$\dot{\phi}(t) \in K^a(\phi(t)); \;\;\phi (0) = x; \;\;
t \in [0, \infty).$$
Denote the class of all such $\phi(\cc)$ (for a given $x$)
by $H^a(x)$.  It will be shown in  Lemma \ref{untiltoo}  that
$H^a(x)$ is nonempty.
Our modified form of the Lyapunov function ($V^a(\cc)$) is defined as follows.
$$
V^a(x) \Df \sup_{\phi \in H^a(x)}\int_0^{\infty}
\eta(|\phi(t)|) dt, \;\; x \in \R^n .
$$
The main step in the proof  is the following
result.  Once this result is proven, parts (a), (b) and (c) of
Definition \ref{def:lyap} used in the statement of Theorem \ref{lyexsts}
follow immediately via one final modification, which
consists of further smoothing, radial linearization, and restriction to $G$,
in exactly the form of \cite{dupwil}(pages 696-697).
Radial linearity of the function thus obtained holds by construction.
Finally, part (d) of Definition \ref{def:lyap} follows immediately from
radial linearity and the fact that the function is $C^2$ on $G\setminus
\{0\}$.

\bt
\label{interm1}
There exist $ a_0 \in (0, \infty)$ such that the following hold
for all $a \in (0,
a_0)$.
\ben
\item  There exists $r \in (0, 1)$, not depending
on $a$ such that $V^a(x) = 0$
for all $x \in B(0,r)$.
\item
 $V^a(\cc)$ is locally Lipschitz on $\R^n$.  In fact, for all $R \in \R^n$
 there exists $\alpha(R) \in (0, \infty)$ and $\ov C(R)$ such that
 for all $x,y \in \R^n$ with $|x| \le R$ and $|x-y| \le \ov C(R)$,
 $$|V^a(x) - V^a(y)| \le \alpha(R) |x-y|, \;\; \forall a \in (0,a_0).$$
 \item
 For
a.e. $x \in \R^n$; $|x| \ge 2$,
$$
\max_{u \in K^a(x)}\inn{DV^a(x)}{u} \le -1.
$$
\item
There exists $D \in (0, \infty)$ such that for all $x \in \R^n$
with $|x| \ge 2$, $V^a(x) \ge \frac{|x| -2}{D}$.
\item
There exists $M \in (0, \infty)$ such that
$$\essinf_{x \in \R^n: |x| \ge M}
\inn{DV^a(x)}{\frac{x}{|x|}} > 0 .$$
\een
\et
In the remaining part of this section we will prove the above
theorem.
The main idea in the proof is that the stability properties
of the trajectories introduced in Proposition \ref{staba}
imply similar properties for the solutions of the differential
inclusion $\dot{\phi}(t) \in K^a(\phi(t)); \;
\phi(0) = x$ for small enough value of $a$.
More precisely, the following result will be shown.
\begin{prop}
\label{prop22-24}
There exist $a_0, T \in (0, \infty)$
such that the following hold.
\ben
\item
Whenever $g(\cc)$ is an absolutely continuous function
such that for some $a \in (0, a_0)$,
$$
\dot{g}(t) \in K^a(g(t))
\quad {\rm a.e.}\; t \in [0, \infty); \quad |g(0)| \le 2^m,
$$
we have that $g(t) = 0$ for all
$t \ge 2^{m+1}T$.
\item  There exist $r \in (0,1)$
such that
whenever $\phi(\cc)$ is an absolutely continuous function and
$$
\dot{\phi}(t) \in K^a(\phi(t)), \quad {\rm a.e.}\; t;\quad |\phi(0)| \le r;
\quad a \in (0, a_0),
$$
we have that $\sup_{0 \le t < \infty}|\phi(t)| \le 1$ .
\een
\end{prop}
We now give the proof of Theorem \ref{interm1}
assuming that Proposition \ref{prop22-24} holds.

A straightforward calculation shows that for a fixed $a \in
(0,\infty)$,
 and $m,M$  positive finite numbers, there exists
a constant $C(a,m,M) < \infty$ such that
\be
\max_{i \in \{0, \ldots, N\}}
\sup_{\beta \in \calC_1}
\;\;\sup_{x,y: m \le |x|, |y| \le M}
|v_i^a(\beta, x) - v_i^a(\beta, y)|
\le C(a,m,M) |x-y|.
\label{vel-lip}
\ee
Furthermore , it is easy to see that there exists $D \in (0, \infty)$
such that
\be
\max_{i \in \{0, \ldots, N\}}
\sup_{\beta \in \calC_1}
\sup_{x \in \R^n}\;\; \sup_{a \in (0, \infty)}|v_i^a(\beta, x)| \le D.
\label{vel-bd}
\ee
\noi{\bf Proof of Theorem \ref{interm1}.}
Let $a_0$ and $r$ be as in Proposition \ref{prop22-24}.
The choice of $r$ implies that if $|x| \le r$
and if $\phi(\cc) \in H^a(x)$ then
$\phi(t) \in B(0,1)$ for all $t \in [0, \infty)$.
This implies that $\eta(|\phi(t)|) = 0$ for all $t$, thus
for such $x$, $V^a(x) = 0$.  This proves part 1.
Now we show the local Lipschitz property in 2.
Let $x \in \R^n$ be such that $|x| \le R$.
Without loss of generality we can assume that
$|x| \ge \frac{r}{2}$ for else local Lipschitz property holds trivially.
From Proposition \ref{prop22-24}(1) it follows that
we can choose $T_0 < \infty$  such that for any $\phi \in H^a(y)$; $|y|
\le R +1$
we have that $\phi(t) = 0$ for all $t \ge T_0$.
For an absolutely continuous trajectory
$\phi: [0, \infty) \mapsto \R^n$, define
 $$\tau^*(\phi) \Df \inf \{ t \in (0, \infty) : \phi(t) \in
B(0,{r/2})\}.$$
Now let $\phi \in H^a(x)$ and $x \in B(0,R)$ be such that
$$
V^a(x) \le \int_0^{\tau^*(\phi)} \eta(|\phi(t)|) dt + \ep.
$$
Note that we could replace $\infty$ by
$\tau^*(\phi)$ in the upper limit of the integral on the right,
because of Proposition \ref{prop22-24}(2).
Let $y \in \R^n$ be such that,
$|y| \le R+1$.

It will be shown in Lemma
\ref{selectb} that
there exist measurable functions
$q_i:[0,\infty) \to [0,1]$; $i = 0, \ldots, N$ and
$\beta: [0,\infty) \to \C_1$ such that $\phi(\cc)$ solves
\beq
\dot{\phi}(t) &=& \sum_{i=0}^N q_i(t) v_i^{a}(\beta(t),\phi(t)),
\;\; {\rm a.e.}\; t \in [0,\infty)\\
\phi(0) &=& x.
\eeq
Now let $\psi(\cc)$ be an absolutely continuous function such
that
for a.e.  $t \in [0, \infty)$
\beq
\dot{\psi}(t) &=& \sum_{i=0}^N q_i(t) v_i^{a}(\beta(t),\psi(t)),
\\
\psi(0) &=& y.
\eeq
Existence of such a $\psi(\cc)$ will be proved in
Lemma \ref{untiltoo}.  Since $\psi \in H^a(y)$, we have that
 $\tau^*(\psi) \le T_0$.
We now claim that if $y$ is sufficiently close to $x$
then both $\phi(\tau^*(\phi)\wedge \tau^*(\psi))$
and $\psi(\tau^*(\phi)\wedge \tau^*(\psi))$ are in $B(0,r)$.
To see this note that as a consequence of
(\ref{vel-lip}) and (\ref{vel-bd}), for $t \in [0, \tau^*(\phi)\wedge
\tau^*(\psi)]$,
$$
|\phi(t) - \psi(t)| \le |y-x| +
C^* \int_0^t
|\phi(s) - \psi(s)| ds,
$$
where $C^* \Df C(a,  \frac{r}{2}, R+1 + T_0D)$.
By an application of Gronwall's inequality we see now that
if $|y-x| \le \frac{r}{2}\exp(-C^*T_0) \equiv \ov C$
then
 \be
\label{anlip}
|\phi(t) - \psi(t)|  \le
\exp(C^*T_0) |y-x| \le
\frac{r}{2}
\ee
 for all $t \in
[0, \tau^*(\phi)\wedge
\tau^*(\psi)].$

This means that for such $y$ both $\phi(\cc)$ and $\psi(\cc)$
are in $B(0,r)$ at time $\tau^*(\phi)\wedge
\tau^*(\psi)$.
Henceforth we will only consider such $y$ (i.e.\ $|y-x| \le \ov
C$).
Note next that
\beq
V^a(x) - V^a(y)
& \le &
\int _0^{\tau^*(\phi) \wedge \tau^*(\psi)}(
\eta(|\phi(t)|) - \eta(|\psi(t)|) dt + \ep.\\
& \le &
\eta_{lip} T_0e^{C^*T_0} |y-x| + \ep,
\eeq
where $\eta_{lip}$
is the Lipschitz constant for
$\eta(|\cc|)$.
Sending $\ep \to 0$ and using the symmetry of the above
calculation
we have that
\be
|V^a(x) - V^a(y)| \le \eta_{lip}e^{C^*T_0} |y-x|
\label{lipV}
\ee
for all $|x| \le R$ and $|y-x| \le \ov C$.
Since $R> 0$ is arbitrary, this proves part 2.

To prove part 3, we will show that
at all points $x$ at which $V^a(\cc)$ is differentiable
and $|x| \ge r$
\be
\label{show}
\max_{u \in K^a(x)} \inn{DV^a(x)}{u} \le -\eta(|x|).
\ee
Fix $R \in [2, \infty)$.
Now let $r \le |x| \le R-1$ and $u \in K^a(x)$.  Then there exist
$q_i \in [0,1]$; $i = 0, \ldots, N$ satisfying
$\sum_{i=0}^N q_i =1$ and $\beta \in \C_1$
such that
$u = \sum_{i=0}^N q_i v_i^{a}(\beta,x).$
Define for $y$ such that $|x-y| < \frac{r}{2}$,
$u(y) \Df \sum_{i=0}^N q_i v_i^{a}(\beta,y).$
In view of (\ref{vel-lip}) there exists $\ti C \equiv C(a, r/2, R+1)$
 such that
$$
|u(x) -u(y)| \le \ti C|x-y|.
$$
Now for a given $y$ such that $|y| \ge \frac{r}{2}$ and
$|x-y| < r/2$.
define $\phi_y(\cc)$ to be the absolutely continuous function
which satisfies
$$
\dot{\phi_y}(t) \in K^a(\phi_y(t)), \;\; \phi_y(0) = y$$
for $t \in [0, \infty)$ and is $\ep$-optimal, i.e.
$$
V^a(y) \le \int_0^{\infty}\eta(|\phi_y(s)|) ds + \ep.
$$
Let $\phi(\cc)$ be an absolutely continuous function such that
$\phi(\cc)$ solves:
$$
\dot{\phi}(t) = u(\phi(t)) ; \;\; \phi(0) = x,$$
$t \in [0, \infty)$.  The existence of such a $\phi$
is again assured from Lemma \ref{untiltoo}.
Now let $\tau_0 > 0$ be such that
for all $t \in [0, \tau_0]$,
$|\phi(t) - x| < \frac{\ov C}{2}$
and
$\tau_0|u| < \frac{\ov C}{2}$.
Now set $y \equiv \phi(\tau_0)$.
Note that since $\ov C < r$, we have that  $|x-y| \le \frac{r}{2}$
and $|y| \ge r/2$.
Consider the following modification of the trajectory
$\phi(\cc)$.
\beq
\ti \phi(t) &=& \phi(t); \;\; t \in [0, \tau_0]\\
\ti \phi(t) & =& \phi_y(t-\tau_0); \;\; t \ge \tau_0
\eeq
Note that by construction $\ti \phi(\cc)$ solves the differential
inclusion:
\beq
\dot{\ti \phi}(t) \in K^a(\ti \phi(t)); \;\; \forall
t \in [0, \infty).
\eeq
Now an argument, exactly as on pages 694-695 of \cite{dupwil}
shows that
$$\frac{V^a(x+ \tau_0u)- V^a(x)}{\tau_0}
\le -\eta(|x|) + O(\tau_0) - \frac{1}{\tau_0} \int_0^{\tau_0}
(\eta(|\ti \phi(s)|) - \eta(|x|)) ds.$$
Taking limit as $\tau_0 \to 0$ we have part 3.

Now we consider part 4.
Let $\phi \in H^a(x)$ and let $\ti \tau \Df \inf \{t: \phi(t) \in
B(0,2)\}$.
Then
we have that
\beq
2 & \ge& |\phi(0)| - |\int_0^{\ti \tau} \dot{\phi}(s) ds| \\
& \ge & |x| - \ti \tau D.
\eeq
Thus $\ti \tau \ge \frac{|x| -2}{D}$.
Since $\eta(|x|) = 1$ for $|x| \ge 2$ we have that
\be
\label{vabd}
V^a(x)  \ge
\int_0^{\ti \tau} \eta(|\phi(s)|) ds
\ge \frac{|x| -2}{D}.\ee
This proves part 4.

Finally, we consider part 5.
We will show that there exists $\alpha \in (0, \infty)$ such that
for all $x \in \R^n$
for which $V^a(x)$ is differentiable, we have that
$$
\inn{DV^a(x)}{\frac{x}{|x|}} \ge \frac{1}{\alpha}(1 -
\frac{2}{|x|}).
$$
This will clearly yield part 5.
Without loss of generality assume that
$|x| \ge \frac{r}{2}$ since otherwise the inequality holds trivially.
In order to show the inequality it suffices to show, in view of part 4, that
\be
\label{mainbd}
\inn{DV^a(x)}{\frac{x}{|x|}} \ge \frac{V^a(x)}{|x|}.
\ee
Now the proof of (\ref{mainbd}) is identical to the proof
of Proposition 3.7 of \cite{dupwil} on observing
that if $\phi \in H^a(x)$ then for $c \in (0, \infty)$,
the trajectory
$\theta^c(\cc)$, defined as
$\theta^c(t) \Df (1+c) \phi(\frac{t}{1+c}),\; t > 0$,
is in $H^a((1+c)x)$.  We omit the details.
\ink

The proof of Theorem \ref{interm1} used in addition to
Proposition \ref{prop22-24}, the following two lemmas.
The first lemma is a classical existence and uniqueness
result, a sketch of whose proof is provided in the appendix,
while the second is a result on measurable selections.
\begin{lemma}
\label{untiltoo}
Let $a_0$ be as in Proposition \ref{prop22-24} and $a \in (0, a_0)$ be fixed.
Let $q_i(\cc)$; $i = 0, \ldots, N$ be measurable
functions from $[0,\infty) \to [0,1]$ such that
$\sum_{i=0}^N q_i(t) = 1$ for all $t \in [0,\infty)$
and $\beta(\cc)$ be a measurable
function from $[0,\infty) \to \C_1$.
  Let $y \in
\R^n$ be arbitrary.
Then there
exists an absolutely continuous function $\phi(\cc)$ on $[0,\infty)$
such that
\beqn
\dot{\phi}(t) &=& \sum_{i=0}^N q_i(t)v_i^{a}(\beta(t),\phi(t));
\;\; {\rm a.e.}\; t \in [0, \infty)\non \\
\phi(0) &=& y.
\label{until}
\eeqn
Furthermore if $\psi(\cc)$ is another absolutely continuous
function solving (\ref{until}), then $\phi = \psi$.
\end{lemma}
\begin{lemma}
Let $a>0$ be fixed and
$\phi(t)$ be an absolutely continuous function on $[0,T]$
such that
$$\dot{\phi}(t) \in K^a(\phi(t)), \;\;\; {\rm a.e.}\; t \in [0,T].$$
Then there exist measurable functions
$q_i:[0,T] \to [0,1]$; $i = 0, \ldots, N$ and
$\beta: [0,T] \to \C_1$ such that $\sum_i q_i(t) = 1$ and
\be
\label{solveb}
\dot{\phi}(t) = \sum_{i=0}^N q_i(t) v_i^{a}(\beta(t),\phi(t)),
\;\; {\rm a.e.} \;t \in [0,T].
\ee
\label{selectb}
\end{lemma}

\noi{\bf Proof:}
Let  $B$ be the subset of $\R^n \times \R^n\setminus \{0\}$
defined as
$$
\{(u,x) \in \R^n \times \R^n\setminus \{0\}:
u = \sum_{i=0}^N q_i v_i^{a}(\beta,x); \;\;
q_i \in [0,1]; i = 0, \ldots, N, \;\; \beta \in \C_1, \;\;
\sum_{i=0}^N
q_i = 1 \}.
$$
Let $\B^N$ be the Borel $\sigma$-field on $[0,1]^{N+1}\times \C_1$.
Define
$
F: B \mapsto \B^N$
 as
$$
F(u,x) = \{(q, \beta): q \in [0,1]^{N+1}; \beta \in \C_1;
\sum_{i = 0}^Nq_i v_i^{a}(\beta,x) = u; \;\; \sum_{i=0}^N q_i =
1\}.
$$
Note that the map $(x, \beta) \to v_i^{a}(\beta,x)$ is continuous
on
$\R^n\setminus \{0\} \times \C_1$ for all $i = 0, 1, \ldots, N$.
This implies that if we have a sequence
$(q_k,\beta_k,u_k,x_k) \to (q,\beta,u,x) \in [0,1]^{N+1} \times
\C_1
\times \R^n \times \R^n\setminus \{0\}$
and
$(q_k,\beta_k) \in F(u_k,x_k)$ for all $k$ then
$(q,\beta)  \in F(u,x)$.
Thus in view of Corollary 10.3, Appendix of \cite{EK}
there exists a measurable selection for $F$, i.e.\ there exists a
measurable map:
$$f: B \to [0,1]^{N+1}\times \C_1$$
such that
$f(u,x) \in F(u,x)$ for all $(u,x) \in B$.
Choose an arbitrary element $(\hat q, \hat \beta) \in\; [0,1]^{N+1}
\times \C_1$ such that $\sum_{i=0}^N q_i =
1$ and extend $f$ to $\hat B \Df B \cup ({\bf 0},0)$ by
setting
$f({\bf 0}, 0) \Df (\hat q, \hat \beta)$.
Now write $f(\cc)$ as $(\{f_i(\cc)\}_{i=0}^N, f_*(\cc))$,
i.e. we denote the first $N+1$ coordinates of the vector function
$f$ by $f_i$; $i = 0, \ldots, N$, and we denote the $N+2$'th coordinate
by $f_*$.
Define
$$
q_i(t) \Df f_i(\dot{\phi}(t), \phi(t)); i = 0,1, \ldots, N; \;\;
{\rm a.e.} \; t \in [0,\infty)$$
$$
\beta(t) \Df f_*(\dot{\phi}(t), \phi(t)); {\rm a.e.}\; t \in [0,\infty).$$
Clearly $q(\cc)$ and $\beta(\cc)$ are measurable functions and by
construction (\ref{solveb}) holds.
\ink

Now we turn to the proof of Proposition \ref{prop22-24}.
The key idea is to relate  the  solutions of the differential
inclusion $\dot{\phi}(t) \in K^a(\phi(t)); \;
\phi(0) = x$, for small enough value of $a$, with
the solutions of the SP
for trajectories with velocity in  $\C_1$.
 The following two results are central in that respect.
 Define for $x \in \R^n$
\beq
K(x) &\Df& \{v \in \R^n : \mbox{there exists a sequence}\;\;
(a_k,x_k,v_k)_{k \ge 1} \subset  (0,1]\times \R^n \times \R^n
\\
&&\ {\rm s.t.}\ a_k \to 0; \;
x_k \to x; \;v_k \to v; \;\mbox{and}\;v_k \in K^{a_k}(x_k) \}.
\eeq
We will denote the closure of the convex hull of $K(x)$ by
$\ov K(x)$.
The first result shows that as $a$ approaches $0$ the solutions
of the  differential inclusion converge to a trajectory
which also solves a differential inclusion given in terms
of $\ov K$.  The proof is quite similar to the proof of
Proposition 3.3 of \cite{dupwil}.  We provide a sketch
in the appendix.
\begin{lemma}
\label{odeKbar}
Consider the sequence
$(x_k, a_k, \phi_k(\cc))_{k \ge 1} \subset
\R^n \times (0,1] \times \C([0,\infty); \R^n)$
such that $x_k \to x$; $a_k \to 0$,
and $\phi_k(\cc) \to \phi(\cc)$ (uniformly on compacts).
Suppose further that each $\phi_k$ is absolutely continuous
and solves the differential inclusion:
$$
\dot{\phi}_k(t) \in K^{a_k}(\phi_k(t));\; {\rm a.e.}\; t\in [0, \infty);
\;\;\; \phi_k(0) = x_k.$$
Then
$\phi(\cc)$ is Lipschitz continuous (and thus absolutely
continuous)
and it solves the differential inclusion:
$$
\dot{\phi}(t) \in \ov K(\phi(t));\; {\rm a.e.}\; t\in [0, \infty);\;\;\;
\phi(0) = x.$$
\end{lemma}
The proposition below provides the connection between the
solutions of  the differential
inclusion $\dot{\phi}(t) \in \ov K(\phi(t)); \;
\phi(0) = x$
and certain solutions to  the SP.
Define
\be
\delta_0 \Df \inf_{(\la,i):\la \in \Lambda,\, i \in \la}
\inn{d^{\la}}{n_i}.
\label{delo}
\ee
By Condition \ref{compS}, $\delta_0 > 0$.
\begin{prop}
Let $\phi: [0, \infty) \to \R^n$ be an absolutely continuous
function
which solves the differential inclusion:
$$
\dot{\phi}(t) \in \ov K(\phi(t)), \quad {\rm a.e.}\ t \in [0,
\infty).$$
Then there exists a $\tau \in [0, {|\phi(0)|}/{\delta_0})$,
a strictly increasing, onto function
$\alpha: [0, \infty) \to [0, \infty)$ and a measurable function
$\ov \beta:[0, \infty) \to \C_1$ such that
$$
\psi(t) \Df \phi(\tau + \alpha(t))
= \Gamma\left(\phi(\tau) + \int_0^{\cc}\ov \beta(u)
  du\right)(t), \;\; t \in (0, \infty).$$
\label{lim-skor}
\end{prop}
Before presenting the proof of this proposition, we show how the
proof of Proposition \ref{prop22-24} follows.
For an absolutely continuous trajectory $\theta: [0, \infty) \mapsto
\R^n$, define
\be
\tau(\theta) \Df \inf\{t \in (0, \infty):\theta(t)
= 0\}.\label{taqudef}
\ee

\noi{\bf Proof of Proposition \ref{prop22-24}.}
First we show that there exist $\ti a_0, T \in (0, \infty)$
such that for all $a \in (0, \ti a_0)$ and  absolutely continuous
$\phi(\cc)$
on $[0,
\infty)$ satisfying
$$
\dot{\phi}(t) \in K^a(\phi); \;\; |\phi(0)| \le 1;
 \;\; {\rm a.e.} \; t \in [0, \infty),
$$
we have that
\be
\label{prelim}
\inf_{0 \le t < T}|\phi(t)| \le1/2.
\ee
We argue by contradiction.
Suppose that there exists a sequences $\{T_k\}_{k \ge 1}$
increasing to $\infty$, $\{a_k\}_{k \ge 1}$ decreasing to $0$
and $\{\phi_k(\cc)\}_{k \ge 1}$ such that for all $k$,
$\phi_k(\cc)$
is absolutely continuous, for {\rm a.e.} $t \in [0, \infty)$,
$$
\dot{\phi}_k(t) \in K^{a_k}(\phi_k(t)); \;\; |\phi_k(0)| \le 1$$
and
$\inf_{0 \le t < T_k}|\phi_k(t)| > 1/2$.
As in Proposition 3.3(i) of \cite{dupwil}
we have that $\{\phi_k; k \ge 1\}$ is precompact in
$C([0, \infty); \R^n)$.
Assume without loss of generality
that
$\phi_k(\cc)$ converges to $\phi(\cc)$ uniformly on compacts.
Clearly $|\phi(0)| \le 1$ and
\be
|\phi(t)| \ge 1/2, \;\; \mbox{for all}\;\; t
\in [0, \infty).
\label{contra2}
\ee
From Lemma \ref{odeKbar} we have that
$\phi(\cc)$ is absolutely continuous and solves the differential
inclusion:
$$
\dot{\phi(t)}\in \ov K(\phi(t))\quad {\rm a.e.}\ t \in [0, \infty).$$
Therefore from Proposition \ref{lim-skor} we have that there exists $\tau
\in [0,{\delta_0}^{-1}]$, a
strictly increasing, onto function $\alpha: [0, \infty) \to
[0, \infty)$ and a measurable function $\ov \beta: [0, \infty) \to \C_1$
such that
for all $t \ge 0$
$$
\phi(\tau + \alpha(t))
= \Gamma\left(\phi(\tau) + \int_0^{\cc}\ov \beta(u) du\right)(t).$$
Now applying Proposition \ref{staba} we have that
$\lim_{t \to \infty} \phi(\tau + \alpha(t)) = 0.$
This is a contradiction to (\ref{contra2}).
Hence (\ref{prelim}) is proven.

Now let $\psi$ be an
absolutely continuous
function
on $[0,
\infty)$ satisfying
$$
\dot{\psi}(t) \in K^a(\psi); \;\; |\psi(0)| \le 2^k;
 \;\; {\rm a.e.} \; t\in [0, \infty).
$$
Assume without loss of generality that $\psi(0) \neq 0$ and
define $\phi(t) \Df 2^{-k} \psi(2^kt)$.
Since $K^a(x) = K^a(\alpha x)$ for all $\alpha > 0$ we have that
$$
\dot{\phi}(t) \in K^a(\phi); \;\; |\phi(0)| \le 1;
 \;\; {\rm a.e.} \; t \in [0, \infty).
$$
Thus
\beq
\inf_{0 \le t \le T2^k}|\psi(t)|
& =& \inf_{0 \le t \le T2^k} 2^k|\phi({t}2^{-k})|\\
& =& 2^k \inf_{0 \le t \le T} |\phi(t)| \\
& \le  & 2^k/2.
\eeq

Now let $g(\cc)$ be as in the proposition.  Then
letting $k = m, (m-1), \ldots, 0, -1, \ldots$
we
have that $\inf_{0 \le t < 2^{m+1}T}|g(t)| = 0$
Since $K^a(0) = {\bf 0}$, we have part 1.

Now we consider part 2.
Let $(\ti a_0, T)$ be as above.
We will show that there exists $a_0 \le \ti a_0$ and $r \in (0,1)$
such that the statement 2 in the proposition holds.
We will once more argue by contradiction.
Suppose that there exist sequences $(a_k, r_k, \phi^k(\cc))$
such that
$$
\dot{\phi^k}(t) \in K^{a_k}(\phi^k(t)); \;\; {\rm a.e.} \; t \in [0,
\infty),
$$
$|\phi^k(0)| \le r_k$, $r_k \to 0$, $a_k \to 0$ and
\be
\label{contradicts}
|\phi^k(t_k)| > c \;\mbox{for some} \;\;t_k \in [0, \infty).
\ee
Assume without loss of generality that $r_k \le 1$ and $a_k \le
a_0$
for all $k \ge 1$.
From 1 we know that
$\tau(\phi^k) \le 2T$ for all $k \ge 1$.
Also note that from the uniform Lipschitz property of
$(\phi^k(\cc))$ and noting that $\phi^k(0) \to 0$
as $k \to \infty$ we have that $T^* \Df \inf_k \tau(\phi^k) > 0$,
since otherwise $\phi^k(t_k)$ converges to $0$ along some
subsequence, which contradicts (\ref{contradicts}).
So now assume without loss of generality that
$\tau(\phi^k) \to T^*$, $t_k \to t^*$ and
$\phi^k(\cc) \to \phi(\cc)$ uniformly on $[0, T^*]$ as $k \to
\infty$.
From Lemma \ref{odeKbar} we have that $\phi(\cc)$ is absolutely
continuous on $[0, T^*]$ and solves the differential inclusion
$
\dot{\phi}(t) \in \ov K(\phi(t)), \;\; {\rm a.e.}\; t \in [0, T^*], \;\;
\phi(0) = 0.$
From Proposition \ref{lim-skor} and Proposition \ref{staba} we then have that
$\phi(t) = 0$ for all $t \in [0, T^*]$.  But on the other hand,
since $\phi^k(t_k) > c$, we have that $\phi(t^*) \ge c$, which is a
contradiction.
This proves part 2 and hence the proposition.
\ink

We now prove Proposition \ref{lim-skor}.
We will need the following three lemmas. The first
lemma characterizes the set $\ov K(x)$ and its proof is similar
to  the proof of Proposition 3.2 of \cite{dupwil} (and is thus
omitted). The second lemma says that a solution of the
differential
inclusion $\dot{\phi}(t) \in \ov K(\phi(t))$,
enters $G$ after some finite time, and then stays within $G$.
The third lemma gives
a representation for a solution to the above differential inclusion.
\begin{lemma}
\label{charK}
For $x \in \R^n \setminus \{0\}$.
\beqn
\ov K(x) & \subset & \conv \{\C_1 \cup \{d_i: i \in \In(x) \},
\;\;\; \mbox{if}\;x \in \pa G \non \\
& \subset & \C_1, \;\; \mbox{if}\;x \in G^0 \non \\
& \subset & \conv\{d^{\lambda}: \lambda \supset \{i: \inn{x}{n_i}
< 0\}\}\;\; \mbox{if}\; x \in G^c. \non \\
\label{char}
\eeqn
\end{lemma}
\begin{lemma}
\label{stayin}
Let $\phi: [0, \infty) \to \R^n$ be an absolutely continuous
function
which solves the differential inclusion:
$$
\dot{\phi}(t) \in \ov K(\phi(t)), \mbox{a.e.}\;\; t \in [0,
\infty).$$
Then the following hold.
\begin{enumerate}
\item
Let $t \in [0, \infty)$ be so that
 $\phi(\cc)$ is differentiable at $t$
and $\phi(t) \not \in G$, then
$$
\frac{d}{dt}
\left[ \min_{i \in \{1, \ldots, N\}}
\inn{\phi(t)}{e_i} \right] \ge \delta_0.$$
\item
If $\phi(0) \in G$ then $\phi(t) \in G$ for all $t \in [0,
\infty)$.
\item
If $\phi(0) \not \in G$ then there exists  $\tau \le
\frac{|\phi(0)|}{\delta_0}$
such that $\phi(t) \in G$ for all $t \ge \tau$.
\end{enumerate}
\end{lemma}

\noi{\bf Proof:}
Parts 2 and 3 follow immediately once 1 is proven.
We now present the proof of (1).
Fix $y \in G^c$.  Define
$\lambda_1(y) \Df \{i: \inn{y}{n_i} < 0\}.$
Note that whenever $\lambda \supset \lambda_1(y)$, we have from
(\ref{delo}) that
$
\inn{d^{\la}}{n_i} \ge \del_0,
\;\;\; \forall i \in \la_1(y).$
This yields the implication:
\be
v \in \conv\{d^{\la}: \la \supset \la_1(y)\}
\;\;\;\rimply \;\;\; \inn{v}{n_i}\ge \del_0, \;\; \forall i \in
\la_1(y).
\label{imply}
\ee
Define
$$\ov \la(y) \Df \{i\in \la_1(y): \inn{y}{n_i} = \min_i\{\inn{y}{n_i}\}\}.$$
Now let $t \in [0,\infty)$ be such that $\phi(\cc)$ is
differentiable
at $t$ and $\phi(t) \in G^c$.
Then by continuity of $\phi(\cc)$ we can choose $\ep > 0$ such
that
for all $0 \le h < \ep$:
$$
\min_{i \in \{1, \ldots, N\}}
\inn{\phi(t+h)}{n_i} = \min_{i \in \ov
  \la(\phi(t))}\inn{\phi(t+h)}{n_i}.
$$
and
\be
\la_1(\phi(t+h)) \supset \la_1(\phi(t)).
\label{contain}
\ee
Next observe that
\beq
\frac{d}{dt}
\left[ \min_{i \in \{1, \ldots, N\}}
\inn{\phi(t)}{n_i} \right]
&=&
\lim_{h\to 0}
\frac{1}{h}\left(\min_{i \in \{1, \ldots, N\}}
\inn{\phi(t+h)}{n_i} - \min_{i \in \{1, \ldots, N\}}
\inn{\phi(t)}{n_i} \right) \\
&=&
\lim_{h\to 0}\frac{1}{h}\min_{i \in  \ov
  \la(\phi(t))}
\left(\inn{\phi(t+h)}{n_i} -
\inn{\phi(t)}{n_i} \right) \\
&=&
\lim_{h\to 0}\min_{i \in  \ov
  \la(\phi(t))}\frac{1}{h}\int_0^h \inn{\dot{\phi}(t+s)}
{n_i}ds \\
&\ge & \del_0,
\eeq
where the last step follows from (\ref{imply})
on observing that in view of Lemma
\ref{charK} and (\ref{contain})
for a.e. $s \in [0,h]$
\beq
\dot{\phi}(t+s)& \in & \ov K(\phi(t+s))\\
& \subset & \conv \{d^{\la}: \la \supset \la_1(\phi(t+s))\}\\
& \subset &  \conv \{d^{\la}: \la \supset \la_1(\phi(t))\}.
\eeq
This proves the lemma.
\ink

The following lemma once more uses a result on measurable selections.
The proof is quite similar to Lemma \ref{selectb}, and a sketch is
given in the appendix.
\begin{lemma}
Let $\phi: [0, \infty) \to \R^n$ be an absolutely continuous
function
such that $\phi(0) \in G$ and $\phi$ solves the differential inclusion:
$$\dot{\phi}(t) \in \ov K(\phi(t)); \;\;\; {\rm a.e.}\; t.$$
Then there exist measurable functions $q_i: [0, \infty) \to
[0,1]$;
$i = 0, 1, \ldots, N$, satisfying the equality
$\sum_{i=0}^N q_i(t) = 1$, and measurable map
$\beta_0:[0,\infty) \to \C_1$ such that
for {\rm a.e.} $t \in [0, \infty)$
$$\dot{\phi}(t) = \sum_{i \in \In(\phi(t))}q_i(t) d_i +
q_0(t) \beta_0(t).$$
\label{select1}
\end{lemma}

\noi{\bf Proof of Proposition \ref{lim-skor}:}
From Lemma \ref{stayin} we know that
$\phi(t) \in G$ for {\rm a.e.} $t > \tau$.
From Lemma \ref{select1} it follows that there exist measurable functions
$q_i:[0, \infty) \to [0,1]$; $i = 0,1, \ldots, N$
and $\beta_0: [0, \infty) \to \C_1$ such that for a.e. $t \ge \tau$
$$
\dot{\phi}(t) = \sum_{i \in \In(\phi(t))} q_i(t) d_i + q_0(t)
\beta_0(t).$$
Let $\{n^{\la}\}_{\la \in \Lambda}$ be as in Remark \ref{dual}.
Define
$\del_{*} \Df \inf_{\la \in \Lambda; i \in \la}
\inn{n^{\la}}{d_i}.$
Also let
$\gamma \Df \sup_{\beta \in \C_1} |\beta|.$
We now claim that for a.e. $t \in [\tau , \infty)$
\be
q_0(t) \ge \frac{\del_{*}}{\del_* + \gamma}.
\label{claim}
\ee
Let $\la \in \Lambda \setminus \{0\}$ be arbitrary.  Define
$$F^{\la} \Df \{x \in \R^n: \In(x) \supset \la\}.$$
Since $\phi(\cc)$ is absolutely continuous and $F^{\la}$ is a
linear subspace of $\R^n$ we have that for a.e. $t$
whenever
$\phi(t)
\in F^{\la}$ we have that $\dot{\phi}(t) \in F^{\la}$.
Thus for a.e. $t$,
$
I_{\{\phi(t) \in F^{\la}\}}\inn{\dot{\phi}(t)}{n^\la} = 0.
$
Now observe that for a.e $t \ge \tau$ such that $\phi(t) \in F^{\la}$:
\beq
0 &=& \inn{n^{\la}}{\dot{\phi}(t)} \\
&= &
n^{\la}\cdot(\dot{\phi}(t) - q_0(t)\beta_0(t)) +
\inn{n^{\la}}{q_0(t)\beta_0(t)}
\\
&=& \inn{n^{\la}}{\sum_{i \in \la}q_i(t)d_i}
+ q_0(t)\inn{n^{\la}}{\beta_0(t)} \\
&\ge&
\del_{*}\sum_{i \in \la}q_i(t) -\gamma q_0(t).
\eeq
This proves (\ref{claim}) for
a.e. $t \ge \tau$ such that $\phi(t) \in F^{\la}$.
Also the claim holds trivially if $\phi(t) \in G^0$ since then
$q_0(t) \equiv 1$.
Now letting $\la$ run over all the subsets of $\Lambda$ we have the claim.
Next define the strictly increasing function
$a: [0, \infty) \to [0, \infty)$ as
$$a(t) \Df \int_0^t q_0(\tau + s) ds; \;\;\; t \in [0, \infty).$$
Also set $\alpha(t) \Df a^{-1}(t)$.
Finally we show that
$\psi(\cc) \Df \phi (\tau + \alpha(\cc))$ solves the SP
for
$$
x(\cc) \Df x(0) + \int_0^{\cc} \ov \beta(s) ds,$$
where $x(0) = \phi(\tau)$ and
$\ov \beta(t) \Df \beta_0(\tau + \alpha(t))$, $t \in [0, \infty)$.
To see this we only need to observe that for a.e. $t \ge 0$
\beq
\dot{\psi}(t) & =& \frac{\dot{\phi}(\tau + \alpha(t))}{q_0(\tau +
  \alpha(t))}
\\
& = & \beta_0(\tau + \alpha(t)) + \sum_{i \in \In(\psi(t))}
\frac{q_i(\tau + \alpha(t))}{q_0(\tau + \alpha(t))}d_i.
\eeq
This proves the lemma.
\ink

\section{Appendix}
\noi{\bf Proof of Lemma \ref{untiltoo}:}
The proof is via Picard iteration method.
Define $\phi^{(0)}(\cc) \equiv y$ on $[0,T]$.
For $k \ge 1$, define for $t \in [0,T]$
$$
\phi^{(k)}(t) \Df y + \sum_{i=0}^N\int_0^t q_i(s)
v_i^{a}(\beta(s),\phi^{(k-1)}(s)) ds.
$$
Note that the boundedness of $q_i(\cc)$ and
$v_i^{a}(\beta(\cc),\cc)$
assures that $\phi^{(k)}(\cc)$ is an equicontinuous family (in
fact uniformly Lipschitz continuous) which
is pointwise bounded on $[0,T]$ for all $T < \infty$.
Thus there exists a subsequential (uniform) limit $\phi(\cc)$.
Clearly $\phi(\cc)$ is Lipschitz continuous and thus absolutely
continuous.
Note that the map $(x, \beta) \to v_i^{a}(\beta,x)$ is continuous
on
$\R^n\setminus \{0\} \times \C_1$.  Therefore
 we have that as $k \to \infty$,
$v_i^{a}(\beta(t),\phi^{(k-1)}(t)) \to v_i^{a}(\beta(t),\phi(t))$
for all $t \in [0, \tau(\phi))$, where
$\tau(\phi)$ is as defined in (\ref{taqudef}).
Now a straightforward application of the dominated convergence
theorem shows that $\phi(\cc)$ solves (\ref{until}) on $[0, \tau(\phi))$
and hence, since $v_i^a(\cc, 0) = {\bf 0}$,  on $[0, \infty)$.
Now let $\phi(\cc)$ and $\psi(\cc)$ be two solutions to
(\ref{until}).  We will show that
\be
\label{alb1}
\phi(t) = \psi(t), \forall \; t \in [0, \tau(\phi)\wedge \tau(\psi)).
\ee
Fix $\ep > 0$.  Then there exists $m> 0$ such that
$\min(|\phi(t)|, |\psi(t)|) > m$ on
$[0, \tau(\phi)\wedge \tau(\psi) - \ep)$.
Also let
$M \Df \sup_{0 \le t \le \tau(\phi)} \max\{|\phi(t)|, |\psi(t)|\}.$
Then from (\ref{vel-lip}) we have that
$$
|\phi(t) - \psi(t)| \le C(a, m, M) \int_0^t
|\phi(s) - \psi(s)| ds, $$
for all $t \in [0, \tau(\phi)\wedge \tau(\psi) - \ep)$.
An application of Gronwall's inequality shows that
$\phi$ and $\psi$ are equal on $[0, \tau(\phi)\wedge \tau(\psi) - \ep)$.
Since $\ep > 0$ is arbitrary, we have
(\ref{alb1}).  This also implies that $\tau(\phi) = \tau(\psi)$
and since both trajectories stay at $0$ once they hit $0$,
we have the desired uniqueness on $[0, \infty)$.
\ink

\noi{\bf Proof of Lemma \ref{odeKbar}.}
The Lipschitz continuity of $\phi$ follows immediately on
observing that for $0 \le s \le t < \infty$ and $k \ge 1$
$
|\phi_k(t) - \phi_k(s)| \le D|t-s|.$
We will show that for all
$T \in [0, \infty)$, $\dot{\phi}(t) \in \ov K(\phi(t)); {\rm a.e.}\; t \in
[0,T]$. Fix $T \in [0, \infty)$.
Define a sequence of probability measures on
$\Om_0 \Df \R^n \times \R^n \times [0,T]$ as follows.  For $f \in
C_b(\R^n \times \R^n \times [0,T])$ define
$$
\int_{\Om_0}f(x,y,t) d\mu_k(x,y,t) \Df
\frac{1}{T}\int_0^T f(\phi_k(s), \dot{\phi}_k(s), s) ds.$$
Since
\be
\label{supphi}
\sup_{k \ge 1, s \in [0,T]}|\phi_k(s)| \le \sup_{k \ge
  1}|x_k|
+ DT \Df  \ov C < \infty \ee
and
$|\dot{\phi}_k(s)| \le D$ a.e.\ $s$, we have that
$\{\mu_k\}_{k \ge 1}$ is a tight family of probability measures.
Without loss of generality assume that $\mu_k$ converges weakly
to
$\mu$.
The sequence $\{\mu_k\}$ gives the following useful
representation for
$\{\phi_k\}$:
$$
\phi_k(t) = x_k + \int_0^t \int_{\R^n \times \R^n}
y d\mu_k(x,y,s).$$
Taking limits in the above equality we have
\be
\phi(t) = x + \int_0^t \int_{\R^n \times \R^n}
y d\mu(x,y,s).
\label{muklim}
\ee
Next note that the marginal distribution of $\mu_k$ in the time
variable
is the normalized Lebesgue measure on $[0,T]$ for every $k$
and thus $\mu$ also has the same marginal distribution.
Therefore there exists $\ti \mu(s, \cc)$, a regular conditional probability
distribution, such that
for $f \in  C_b(\R^n \times \R^n \times [0,T])$
\be
\int_{\Om_0} f(x,y,t) d\mu(x,y,t) =
\frac{1}{T}\int_0^T \left(\int_{\R^n \times \R^n} f(x,y,t)\ti \mu(t, dx
  ,dy)\right)
dt.
\label{rcpd}
\ee
Thus from (\ref{muklim}) we have that
$$
\phi(t) = x + \int_0^t \left(\int_{\R^n \times \R^n}
y \ti \mu(s, dx, dy)\right) ds.$$
This shows that for a.e.\ $t \in [0,T]$,
\be
\dot{\phi}(t)
=
\int_{\R^n \times \R^n} y \ti \mu(t, dx, dy).
\label{differen}
\ee
Using the upper-semi continuity of the set $K(x)$, it follows
as in \cite{dupwil} (see pages 687-689) that
the support of
 $\ti \mu(t, dx, dy)$ is contained in
$\{\phi(t)\} \times \ov K(\phi(t))$.  Thus we have from (\ref{differen}),
on noting that
$K(x) \subset  \ov K(x)$ and $\ov K(x)$ is a closed convex set, that
$
\dot{\phi}(t)
\in  \ov K(\phi(t)).
$
This proves the lemma.
\ink

\noi{\bf Proof of Lemma \ref{select1}.}
For $\lambda \in \Lambda$ define
$$
B^{\lambda}
\Df \{u \in \R^n: u = \sum_{i \in \lambda}
q_i d_i + q_0 \beta; \;\;
\sum_{i \in \lambda}q_i + q_0 = 1; \;\; q_i \ge 0; \;\; \beta \in
\C_1\}.
$$
Denote the class of Borel subsets of $[0,1]^{|\lambda| + 1}
\times \C_1$ by
$\B^{|\lambda|}$.
Define the set-valued map $F^{\lambda}: B^{\lambda}
\to \B^{|\lambda|}$ as follows.  For $u \in B^{\lambda}$
$$F^{\lambda}(u)
\Df \{(q,\beta): q \equiv (q_i)_{i\in \lambda \cup \{0\}}
\in [0,1]^{|\lambda| + 1}; \beta \in \C_1; \;\;
\sum_{i \in \lambda}q_i + q_0 = 1; \mbox{and} \sum_{i \in \lambda}
q_i d_i + q_0 \beta = u
\}.$$
We would like to show that there exists a measurable selection for
$F^{\lambda}$,
i.e.\ there exists a measurable map:
$$f^{\lambda}: B^{\lambda} \to [0,1]^{|\lambda| + 1} \times
\C_1$$
such that for all $u \in B^{\lambda}$,
$f^{\lambda}(u) \in F^{\lambda}(u).$
In order to show this it will suffice to show (in view of
Corollary 10.3, Appendix, \cite{EK}) that
if $(q_k, \beta_k) \in F^{\lambda}(u_k)$ and $u_k \to u$ then
the sequence $(q_k,\beta_k)_{k \ge 1}$ has a limit point in
$F^{\lambda}(u)$.
But this is an immediate consequence of the compactness of
$[0,1]^{|\lambda| + 1} \times
\C_1$.
Now fix such a measurable selection for every $\lambda \in
\Lambda$.
Set $$f^{\lambda}(\cc)
\equiv ((f_i^{\lambda}(\cc))_{i \in \lambda \cup \{0\}}, f_{vel}^{\lambda}),
$$
where $f_i^{\lambda}: B^{\lambda} \to [0,1]$ for $i \in \lambda \cup \{0\}$
and $f^{\lambda}_{vel}: B^{\lambda} \to \clc_1$ are the coordinate maps
defined in the obvious way.
Let $\phi(\cc)$ be as in the statement of the theorem.
Define for all $t$ for which $\In(\phi(t)) = \lambda$ and
$\dot{\phi}(t) \in \ov K(\phi(t))$,
$$
q_0(t)\Df f_0^{\lambda}(\dot{\phi}(t));\;
q_i(t) \Df f_i^{\lambda}(\dot{\phi}(t));\; i \in \la;\;
q_i(t) \Df 0;\;\; i \not \in \la \;
\mbox{and}\;\beta_0(t) \Df f_{vel}^{\lambda}(\dot{\phi}(t))\;.
$$
Thus letting $\lambda$ vary over all the subsets of $\Lambda$ we
have
a.e. defined measurable functions $(q_i(\cc))_{i = 0, 1, \ldots,
N}$, $\beta(\cc)$ as required in the statement of the lemma.
\ink

\newpage

\bibliographystyle{plain}

\end{document}